\documentclass[leqno,final]{siamltex}
\usepackage{amsmath}
\usepackage{graphicx}
\usepackage{amssymb}
\usepackage{color}

\newtheorem{remark}{Remark}

\newcommand{\cT}{\mathcal{T}}
\renewcommand{\div}{\mbox{\rm div\,}}

\newcommand{\cS}{\mathcal{S}}
\newcommand{\cK}{\mathcal{K}}

\newcommand{\cV}{\mathcal{V}}
\newcommand{\cW}{\mathcal{W}}

\newcommand{\mP}{\mathbb{P}}

\begin{document}

\title{Analysis of Fully Discrete Mixed Finite Element Scheme for Stochastic Navier-Stokes Equations
	with Multiplicative Noise}
\markboth{HAILONG QIU}{MIXED FEMS FOR STOCHASTIC NAVIER-STOKES EQUATIONS}

\author{
Hailong Qiu\thanks{School of Mathematics and Physics, Yancheng Institute of Technology,
Yancheng, 224051, China. ({\tt qhllf@163.com})
The work of the this author was partially supported by the NSF grant 11701498.}
}

\maketitle

\begin{abstract}
This paper is concerned with stochastic incompressible Navier-Stokes equations with multiplicative noise in two dimensions with respect to periodic boundary conditions.
Based on the Helmholtz decomposition of the multiplicative noise, semi-discrete and fully discrete time-stepping algorithms
are proposed. The convergence rates for mixed finite element methods based time-space approximation
with respect to convergence in probability for the velocity and the pressure are obtained.
 Furthermore, with establishing some stability and using the negative norm technique,
the partial expectations of the $H^1$
and $L^2$ norms of the velocity error are proved to converge optimally.
\end{abstract}



\section{Introduction}\label{sec-1}
In this paper, we consider the following time-dependent stochastic incompressible Navier-Stokes equations:
\begin{subequations}\label{eq1.1}
\begin{alignat}{2} \label{eq1.1a}
du &=\bigl[\nu\Delta u-(u\cdot\nabla)u-\nabla p\bigr] dt +G(u) dW &&\qquad\mbox{a.s. in}\, D_T:=(0,T)\times D,\\
\div u &=0 &&\qquad\mbox{a.s. in}\, D_T,\label{eq1.1b}\\
u(0)&=u_0 &&\qquad\mbox{a.s. in}\, D,\label{eq1.1d}
\end{alignat}
\end{subequations}
 where $T>0$ denotes time, $\nu>0$ is the viscosity of the fluid, ${u}$ and $p$ denote respectively the velocity and the pressure
 of the problem \eqref{eq1.1} which are spatially periodic with period $L>0$,
and $D=(0,L)^2\subset \mathbb{R}^2$ is a period of the periodic
domain with boundary $\partial D$ and $u_0$ denotes a given initial datum.
Here we assume that $\{W(t); t\geq 0\}$ is an $[L^2(D)]^2$-valued
$Q$-Wiener process. The noise is not divergence-free (i.e., $\div G(u)\neq0$). 

The stochastic system \eqref{eq1.1} can take into account
noise term in the sense of physical or numerical uncertainties
and thermodynamical fluctuations. 
 In \cite{BT1973}, Bensoussan and Temam  started to study the stochastic Navier-Stokes
in mathematical investigation.
The paper \cite{Flandoli_Gatarek95} by Flandoli and Gatarek developed a fully stochastic
theory to prove the existence of a martingale
solution.
This paper \cite{H06} investigated the ergodic
properties for the stochastic Navier-Stokes equations with degenerate noise.
 In the last twenty years, there is a large amount of literature
about the analysis of problem \eqref{eq1.1}.
 We refer
  to \cite{Bensoussan95,LRJ03,BBM14,BM18,CP12,D12,Feng_Qiu18,Feng20}
 and the references therein for
 detailed discussions of the stochastic incompressible Navier-Stokes
equations.

The paper \cite{BCP12} by Brze\'{z}niak, {\em et al.} proposed two fully discrete finite element schemes for
 the stochastic Navier-Stokes equations with
multiplicative noise. By using the compactness argument, the authors analyzed the convergence for the velocity field
to weak martingale solutions in 3D and to strong solutions in 2D.
In \cite{CP12}, Carelli and Prohl studied implicit and semi-implicit
fully schemes for the stochastic Navier-Stokes problem.
 The result in \cite{CP12} is convergence
 of rate (amost) $\frac{1}{4}$ in time and linear convergence in space for the velocity.
However, the convergence of the pressure was not given.
In work \cite{BBM14}, the authors proposed an iterative splitting scheme for stochastic Navier-Stokes equations and
 established a strong convergence in probability in the 2D case. In \cite{BM18}, the authors
studied another time-stepping semi-discrete scheme and derived strong $L^2$ convergence for the velocity. In \cite{Hausenblas19},
 Hausenblas and Randrianasolo proposed a time semi-discrete scheme of stochastic 2D Navier-Stokes equations with 
 penalty-projection method. As noted in \cite{Hausenblas19}, the result is convergence
 of rate (amost) $\frac{1}{4}$ in time for the velocity and the pressure.
In paper \cite{Feng_Qiu18}, Feng and Qiu developed a fully discrete mixed finite element scheme of the
time-dependent stochastic Stokes equations with multiplicative noise and established strong
convergence with rates not only for the velocity but also for the pressure.
The paper \cite{Feng20}, by Feng, {\em et al.} proposed a new fully discrete mixed finite element scheme of the
time-dependent stochastic Stokes equations with multiplicative noise and obtained optimal strong
convergence with rates for both the velocity and the pressure.
In a very recent paper \cite{Breit21}, Breit and Dodgson  considered
a fully discrete time-space finite element scheme and proved strong convergence
with rates for the velocity. The result in \cite{Breit21} is convergence
 of rate (amost) $\frac{1}{2}$ in time and linear convergence in space.
 The error estimate of the velocity field $u$ and its
 time-space numerical solution $u_h^n$ reads as: assume that $Lk\leq (-\epsilon\log h)^{-1}$ for some $L>0$, then for any $\alpha<\frac{1}{2}$
\begin{align}
\mathbb{E}\Bigl[\Omega_{k,h}^{\epsilon}\Bigl(\max_{1\leq n\leq N} \|u(t_n)-u^n_h\|^2_{L^2_x}
&+ k\sum_{n=1}^{N}\|\nabla (u(t_n)-u^n_h)\|^2_{L^2_x}\Bigr) \Bigr]\\ \nonumber
& \leq C\bigl(k^{2\alpha-\epsilon}+h^{2}\bigr),
\end{align}
 where $\Omega_{k,h}^{\epsilon}\subset\Omega$ with $\mathbb{P}(\Omega\backslash\Omega_{k,h}^{\epsilon})\rightarrow 0$
 as $k,h\rightarrow 0$.

 The primary goal of this paper is twofold. Our first goal is to develop an optimally convergent fully discrete finite element
  scheme with inf-sup stability. Our main idea, which is partly used in references \cite{Feng_Qiu18,Feng20},
  is to use the Helmholtz decomposition for the driving multiplicative noise at each time step,
   and then solve velocity and pressure.
  We propose new semi-discrete and fully discrete time-stepping algorithms for problem \eqref{eq1.1} and
 prove the convergence of the velocity and the pressure for the fully discrete scheme for the stochastic Navier-Stokes equations.
The second goal is to prove strong optimal $H^1$ convergence first, and then to obtain $L^2$ convergence of the fully discrete scheme
with the negative norm technique.
To the best of our knowledge, it is the first time that strong optimal $L_t^{\infty}-H_x^1/L_x^2$ convergence
 of the discrete solution $(u_h^n, p_h^n)$
to a fully discrete system of the stochastic Navier-Stokes equations
has been established.
 The highlight of this paper (see section \ref{sec-4}) is to derive the error estimates for the numerical
solution as follows: for any $\alpha<\frac{1}{2}$
\begin{align}\label{eq1.2}
\mathbb{E}\Bigl[\textbf{1}_{\Omega_{k,h}^{\epsilon}\cap\Omega_{h,h}^{\epsilon}\cap\Omega_{h}^{\epsilon}\cap\Omega_{\tau}^{\epsilon}}\Bigl( \|\nabla (u(t_n)-u^n_h)\|^2_{L^2_x}
\Bigr)\Bigr]&\leq C\bigl(k^{2\alpha-\epsilon}+h^{2-4\epsilon}\bigr),   \\\label{eq1.333}
\mathbb{E}\Bigl[\textbf{1}_{\Omega_{k}^{\epsilon}\cap\Omega_{h}^{\epsilon}}\Bigl(\Bigl\|\int_0^{t_m} p(s)\,ds -k\sum^m_{n=1}p^n_h \Bigr\|^2_{L^2_x}\Bigr)\Bigr] &\leq C\bigl(k^{2\alpha-\epsilon}+h^{2-2\epsilon}\bigr),\\\label{eq1.3}
\mathbb{E}\Bigl[\textbf{1}_{\Omega_{k,h}^{\epsilon}\cap\Omega_{h,h}^{\epsilon}\cap\Omega_{h}^{\epsilon}\cap\Omega_{k}^{\epsilon}\cap\Omega_{\kappa_0}\cap\Omega_{\kappa}}\Bigl( \|u(t_n)-u^n_h\|^2_{L^2_x}
\Bigr)\Bigr]&\leq C(\kappa_0,\kappa_1)\bigl(k^{2\alpha-4\epsilon}+h^{4-8\epsilon}\bigr).
\end{align}

The remainder of this paper is organized as follows. In Section \ref{sec-2}, we introduce some
function and space notation for problem \eqref{eq1.1} and obtain a few preliminary results.
In Section \ref{sec-3}, we propose the semi-discrete scheme for problem \eqref{eq1.1} and
derive some optimal error estimates for both the velocity and pressure approximations.
In Section \ref{sec-4}, we prove some optimal error estimates for the fully
discrete scheme for problem \eqref{eq1.1} with the negative norm technique.
 In Section \ref{sec-5}, some numerical results are given to validate the theoretical error estimates.
\section{Preliminaries}\label{sec-2}
\subsection{Notation and assumptions}\label{sec-2.1}
Standard function and space notation will be adopted in this paper. Let $H^m(D)$ $(m\geq 0)$
denote the standard Sobolev space, and $\|\cdot\|_{H^m_x}$ denotes its norm.
Let $H^m_{per}(D)$ be the subspace of $H^m(D)$ consisting of $\mathbb{R}^2$-valued periodic function.
Let $(\cdot,\cdot):=(\cdot,\cdot)_D$ denote the standard $L^2$-inner product.
We also let $(\Omega,\mathfrak{F},\mathfrak{F}_t,\mP)$ be a stochastic basis with a complete right continuous filtration. For a given random
 variable $v$ defined on $(\Omega,\mathfrak{F},\mathfrak{F}_t,\mP)$,
let $\mathbb{E}[v]$ denote the expected value of $v$. 
Let $X$ denote a normed vector space $X$ with norm $\|\cdot\|_{X}$.  Define the following Bochner space
\[
L^p(\Omega,X):=\bigl\{v:\Omega\rightarrow X;\, \mathbb{E}[\|v\|_{X}^p] <\infty \bigr\}
\]
and the norm
\[
\|v\|_{L^p(\Omega,X)}:=\Bigl(\mathbb{E}\Bigl[\|v\|_X^p\Bigr]\Bigr)^{\frac1p}, \qquad 1<p<\infty.
\]
We also define some special space notation as follows:
\begin{align*}
\cV:=[H^1_{per}(D)]^2,\quad
\cW:=\{q\in L^2_{per}(D); \, (q,1)_D=0 \},  \quad
\cV_0:=\bigl\{v\in \cV;\,\div v=0 \mbox{ in }D \bigr\}.
\end{align*}

Let $\cK:=L_0(L_{per}^2(D);[L_{per}^2(D)]^2)$ denote the Banach space of linear operators
from $L_{per}^2(D)$ to $[L_{per}^2(D)]^2$ with finite Hilbert-Schmidt norms denoted by $\|\cdot\|_{\cK}$.
As it is noted \cite{PZ92} that the stochastic integral $\int_0^t\varphi(s)dW(s)$ is an $\{\mathfrak{F}_t\}$-martingale and
the following It\^o's isometry holds:
\begin{equation}\label{eq2.3}
\mathbb{E}\Bigl[\Bigl\|\int_0^T\varphi(s)\,dW(s)\Bigr\|_{L^2_x}^2\Bigr]
=\mathbb{E}\Bigl[\int_0^T\|\varphi(s)\|^2_{\cK}\,ds\Bigr].
\end{equation}

In this paper we assume that $G:[0,T]\times [L_{per}^2(D)]^2\rightarrow L^2(\Omega,\cK)$  satisfies the following conditions:
\begin{subequations}\label{eq2.4}
\begin{align}\label{eq2.4a}
\|G(v)-G(w)\|_{L_2(\cK,L^2_x)}&\leq C\|v-w\|_{L^2_x},\ &&\forall  v, w\in [L_{per}^{2}(D)]^2,\\\label{eq2.4b}
\|G(v)\|_{L_2(\cK,H^{i}_x)} &\leq C \bigl(1+\|v\|_{H^{i}_x}\bigr),\  &&\forall v\in [H_{per}^{i}(D)]^2, \, i=1,2,\\\label{eq2.4c}
\|D^{j}G(v)\|_{L_2(\cK,\mathcal{L}([L^2(D)]^j;L^2(D))} &\leq C,\ && \forall v\in [L_{per}^{2}(D)]^2, \, j=1,2.
\end{align}
\end{subequations}

We introduce the Helmholtz projection \cite{Girault_Raviart86} $P_H:[L_{per}^2(D)]^2\rightarrow \cV$ defined by $P_Hv=\eta$ for every $v\in [L_{per}^2(D)]^2$,
where $(\eta,\xi)\in \cV\times [H_{per}^1(D]^2)/\mathbb{R}$ is a unique decomposition
 such that
$$v=\eta+\nabla\xi,$$
and $\xi\in [H_{per}^1(D)]^2/\mathbb{R}$ satisfies the following problem:
\begin{equation}\label{eq2.5}
(\nabla\xi,\nabla q)=(v,\nabla q), \qquad \forall q\in \cV.
\end{equation}

\subsection{Definition of weak solutions}\label{sec-2.2}
In this subsection we first recall the weak solution definition for problem \eqref{eq1.1}, and refer to \cite{Capinski91,Capinski93,Feng_Qiu18,Feng20}.
We then introduce some regularity of the velocity and the pressure.  

\begin{definition}\label{def2.1}
Assume that $(\Omega,\mathfrak{F},(\mathfrak{F}_t)_{t\geq 0},\mathbb{P})$ is a given stochastic basis and
$u_0$ is an $\mathfrak{F}_0$-measurable
random variable. Then $(u,p)$ is called a weak pathwise solution to problem \eqref{eq1.1} if

(i) the velocity and the pressure $(u,p)$ is $\mathfrak{F}_t$-adapted and
 \begin{align*}
 &u\in C([0,T];[L_{per}^2(D)]^2)\cap L^2([0,T;\cV]),&&\quad \mathbb{P} \mbox{-a.s.},\\
&p\in  W^{-1,\infty}(0,T; \cW),&&\quad \mathbb{P}\mbox{-a.s.}.
 \end{align*}

(ii) the problem \eqref{eq1.1} satisfies
\begin{subequations}\label{eq2.6}
\begin{align}\label{eq2.6a}
\bigl(u(t),v \bigr) &+\int_0^t \bigl[a(u(s),v)+ b(u(s),u(s),v)\bigr]\,ds -d\Bigl(v,\int_0^t p(s)\,ds\Bigr) \\ \nonumber
&=(u_0,v)+\Bigl(\int_0^t G(u(s))\, dW(s),v\Bigr) \quad\forall\, v\in \cV,\\\label{eq2.6b}
d(u,q)&=0,\quad\forall\, q\in \cW,
\end{align}
\end{subequations}
 holds $\mathbb{P}$-a.s. for all $t\in (0,T]$. Where the bilinear forms $a(\cdot,\cdot)$ and $d(\cdot,\cdot)$ are defined
 \begin{align*}
a(v,w) &:= \nu\bigl(\nabla v, \nabla w\bigr),&&\ \forall\, v,w\in \cV, \\
d(v,q)&:=\bigl(\nabla\cdot v, q \bigr),&&\  \forall\, v\in \cV,\,  q\in \cW,
 \end{align*}
and the nonlinear form $b(\cdot,\cdot,\cdot)$ is defined as follows:
\begin{align*}
b(w,u,v)&:=\bigr(w\cdot\nabla u, v\bigr)+\frac{1}{2}\bigr((\nabla\cdot w)u, v\bigr), \ \forall\,  w, u, v \in \cV.
\end{align*}
\end{definition}

Using the similar idea in \cite{Feng20}, problem \eqref{eq2.6} can be considered as a mixed formulation for problem \eqref{eq1.1}.
Thus, we introduce a new pressure $r:=p-\xi(u)\, dW$,
where we apply the Helmholtz decomposition $G(u)=\eta(u)+\nabla \xi(u)$, where $\xi(u) \in \cV$, $\mathbb{P}$-a.s. such that
 \begin{align*}
(\nabla\xi(u),\nabla\phi)= (G(u),\nabla\phi), \qquad \forall\, \phi \in \cV.
\end{align*}
By the elliptic regularity \cite{Girault_Raviart86}, we have
 \begin{align}\label{eq2.77}
\|\nabla\xi(u)\|_{L^2_x}&\leq C\|G(u)\|_{L^2_x}, \\\label{eq2.78}
\|\nabla\xi(u)\|_{H^2_x/\mathbb{R}}&\leq C\|\nabla\cdot G(u)\|_{L^2_x}.
\end{align}
\begin{definition}\label{def2.2}
Assume that $(\Omega,\mathfrak{F},(\mathfrak{F}_t)_{t\geq 0},\mathbb{P})$ is a given stochastic basis and
$u_0$ is an $\mathfrak{F}_0$-measurable
random variable. Then $(u,p)$ is called a weak pathwise solution to problem \eqref{eq1.1} if

(i) the velocity and the pressure $(u,r)$ is $\mathfrak{F}_t$-adapted and
 \begin{align*}
 &u\in C([0,T];L_{per}^2(D))\cap L^2([0,T;\cV]),&&\quad \mathbb{P} \mbox{-a.s.},\\
&r\in  W^{-1,\infty}(0,T; \cW),&&\quad \mathbb{P}\mbox{-a.s.}.
 \end{align*}

(ii) the problem \eqref{eq1.1} satisfies
\begin{subequations}\label{eq2.7}
\begin{align}\label{eq2.7a}
\bigl(u(t),v \bigr) &+\int_0^t \bigl[a(u(s),v)+ b(u(s),u(s),v)\bigr]\,ds -d\Bigl(v,\int_0^t r(s)\,ds\Bigr) \\ \nonumber
&=(u_0,v)+\Bigl(\int_0^t \eta(u(s))\, dW(s),v\Bigr) \quad\forall v\in \cV,\\\label{eq2.7b}
d(u,q)&=0,\quad\forall q\in \cW,
\end{align}
\end{subequations}
holds $\mathbb{P}$-a.s. for all $t\in (0,T]$, where $\eta(u(s)):=G(u(s))-\nabla\xi(u(s))$.
\end{definition}

\begin{remark}
By using similar techniques in \cite{LRJ03,Feng20}, we note that $\int^t_{0}pds,\int^t_{0}rds \in L^2(\Omega,$ $L^2(0,T,\cW))$.
\end{remark}

The next Lemma follows from \cite{Breit21}.
\begin{lemma}\label{lem2.1}
(i), Let $u_0\in L^l(\Omega, [L^2_{div}(D)]^2)$ for some $l\geq2$ and let $G$ satisfy
\eqref{eq2.4}. Then there exists a constant $C>0$, such that
\begin{align*}
\mathbb{E}\Bigl[\sup_{0\leq t\leq T}\|u(t)\|^2_{L^2_x} &+\int_0^T\|\nabla u(s)\|_{L^2_x}^2ds\Bigr]^{\frac{l}{2}}
\leq C.
\end{align*}

(ii), Let $u_0\in L^l(\Omega, \cV_0)$ for some $l\geq2$ and let $G$ satisfy
\eqref{eq2.4}. Then there exists a constant $C>0$, such that
\begin{align*}
&\mathbb{E}\Bigl[\sup_{0\leq t\leq T}\|\nabla u(t)\|^2_{L_x^2}+\int_0^{T}\|\nabla^2 u(s)\|_{L^2_x}^2ds\Bigr]^{\frac{l}{2}}
\leq C.
\end{align*}

(iii),  Let $u_0\in L^l(\Omega, [H^{2}(D)]^2\cap\cV_0)\cap L^{5l}(\Omega, \cV_0)$ for some $l\geq2$ and let $G$ satisfy
\eqref{eq2.4}. Then there exists a constant $C>0$, such that
\begin{align*}
&\mathbb{E}\Bigl[\sup_{0\leq t\leq T}\|\nabla^2 u(t)\|^2_{L_x^2}+\int_0^{T}\|\nabla^3 u(s)\|_{L^2_x}^2ds\Bigr]^{\frac{l}{2}} \leq C.
\end{align*}
\end{lemma}

We finish this section by establishing some
regularity of pressure of various spatial norms.
For the reader's convenience, we here give theirs proofs.

\begin{lemma}\label{lem2.2}
(i), Let $u_0\in L^l(\Omega, L^2(D))$ for some $l\geq2$ and let $G$ satisfy
\eqref{eq2.4}. Then there exists a constant $C>0$, such that
\begin{align}\label{eq2.8}
\mathbb{E}\Bigl[\Bigl\|\int_{0}^{t}p(s)ds\Bigl\|^2_{L^2_x}\Bigl]^{\frac{l}{2}}  &\leq C,\\\label{eq2.9}
\mathbb{E}\Bigl[\Bigl\|\int_{0}^{t}\nabla p(s)ds\Bigl\|^2_{L^2_x}\Bigl]^{\frac{l}{2}}  &\leq C.
\end{align}
where $(u,p)$ is the weak pathwise solution to \eqref{eq1.1}, cf. Definition \ref{def2.1}.

(ii), Let $u_0\in L^l(\Omega, L^2(D))$ for some $l\geq2$ and let $G$ satisfy
\eqref{eq2.4}. Then there exists a constant $C>0$, such that
\begin{align}\label{eq2.100}
\mathbb{E}\Bigl[\Bigl\|\int_{0}^{t}r(s)ds\Bigl\|^2_{L^2_x}\Bigl]^{\frac{l}{2}}  &\leq C,\\\label{eq2.101}
\mathbb{E}\Bigl[\Bigl\|\int_{0}^{t}\nabla r(s)ds\Bigl\|^2_{L^2_x}\Bigl]^{\frac{l}{2}}  &\leq C.
\end{align}
where $(u,r)$ is the weak pathwise solution to \eqref{eq1.1}, cf. Definition \ref{def2.2}.
\end{lemma}
\begin{proof}  For \eqref{eq2.8}, from \eqref{eq2.6a}, we have
 \begin{align}\label{eq2.99}
d\Bigl(v, \int_0^{t}p(s)ds\Bigr)&=\bigl(u(t),v\bigr)- \bigl(u_0,v\bigr)+\int_0^{t}\Bigl[a(u(s),v)
\\ \nonumber
&\quad+\bigl(b(u(s),u(s),v) \Bigr]\,ds-\Bigl(\int_0^{t} G(u(s))\, dW(s),v\Bigr).
\end{align}
Using the Young's inequality, the Poincar\'{e} inequality and the H\"{o}lder inequality, one finds that
 \begin{align*}
d\Bigl(v, \int_0^{t}p(s)ds\Bigr)&\leq C\bigl(\|u_{0}\|_{L^2_x} +\|u(t)\|_{L^2_x}\bigr)\|\nabla v\|_{L^2_x}
+C\int_0^{t}\|\nabla u(s)\|_{L^2_x}ds\|\nabla v\|_{L^2_x}\\
&\quad
+C\int_0^{t} \bigl(\|u(s)\|^2_{L^2_x} +\|\nabla u(s)\|^2_{L^2_x}\bigr)ds\|\nabla v\|_{L^2_x}\\
&\quad
+C\Bigl\|\int_0^{t}G(u(s))\, dW(s)\Bigr\|_{L^2_x}\|v\|_{L^2_x}.
\end{align*}
By the well-known inf-sup condition \cite{Girault_Raviart86},  it follows that
 \begin{align*}
\beta\Bigl\|\int_0^{t}p(s)ds\Bigr\|_{L^2_x}&
\leq C\bigl(\|u_{0}\|_{L^2_x} +\|u(t)\|_{L^2_x}\bigr)
+C\int_0^{t}\|\nabla u(s)\|_{L^2_x}ds \\
&\quad
+C\int_0^{t}\bigl( \|u(s)\|^2_{L^2_x} +\|\nabla u(s)\|^2_{L^2_x}\bigr)ds\\
&\quad
 +C\Bigl\|\int_0^{t}G(u(s))\, dW(s)\Bigr\|_{L^2_x}.
\end{align*}
Taking the expectation,  using It\^{o}'s isometry, \eqref{eq2.4b} and Lemma \ref{lem2.1}, which lead to the desired result.

For \eqref{eq2.9}, from \eqref{eq2.6a} and using the H\"{o}lder inequality, the Young's inequality and the Poincar\'{e} inequality, we get
 \begin{align*}
\Bigl(\int_0^{t}\nabla p(s)ds,v\Bigr)&=\Bigl[\bigl(u(t),v\bigr)-\bigl(u_0,v\bigr)+\int_0^{t}\Bigl[a(u(s),v)+ \\ \nonumber
&\quad+b(u(s),u(s),v)\Bigr]\,ds-\Bigl(\int_0^{t} G(u(s))\, dW(s),v\Bigr)\\&
\leq C\bigl(\|u_{0}\|_{L^2_x} +\|u(t)\|_{L^2_x}\bigr)\|v\|_{L^2_x}
+C\int_0^{t}\|\nabla^2 u(s)\|_{L^2_x}ds\|v\|_{L^2_x}\\
&\quad
+C \int_0^{t}\Bigl(\|u(s)\|^2_{L^2_x} +\|\nabla u(s)\|^2_{L^2_x}+\|\nabla^2 u(s)\|^2_{L^2_x}\Bigl)ds\|v\|_{L^2_x}\\
&\quad
+C\Bigl\|\int_0^{t}G(u(s))\, dW(s)\Bigr\|_{L^2_x}\|v\|_{L^2_x}.
\end{align*}
With the definition of $L^2$ norm, we obtain
 \begin{align*}
\Bigl\|\int_0^{t}\nabla p(s)ds\Bigr\|_{L^2_x}&=\sup_{v\neq 0\in L^2_x}\dfrac{\Bigl(\int_0^{t}\nabla p(s)ds,v\Bigr)}{\|v\|_{L^2_x}}\\&
\leq C\Bigl(\|u_{0}\|_{L^2_x} +\|u(t)\|_{L_x}^2\Bigr)
+C\int_0^{t}\|\nabla^2 u(s)\|_{L^2_x}ds \\
&\quad
+C\int_0^{t}\Bigl(\|u(s)\|^2_{L^2_x} +\|\nabla u(s)\|^2_{L^2_x}+\|\nabla^2 u(s)\|^2_{L^2_x}\Bigl)ds\\
&\quad
 +C\Bigl\|\int_0^{t}G(u(s))\, dW(s)\Bigr\|_{L^2_x}.
\end{align*}
Taking the expectation, using It\^{o}'s isometry and \eqref{eq2.4b} and Lemma \ref{lem2.1}, we get the desired result.

 Similarly, using \eqref{eq2.4b} and the definition of $\eta(u)$, \eqref{eq2.77} and Lemma \ref{lem2.1}, the results \eqref{eq2.100}--\eqref{eq2.101} hold.
 The proof is complete.
\end{proof}

\section{Semi-discrete time-stepping scheme}\label{sec-3}
In this section we establish semi-discrete time-stepping scheme for the mixed formulation
\eqref{eq2.7}. Then we analyze the error estimates  for the velocity and the pressure.

Let $N$ be a positive integer and $0=t_0<t_1<\ldots<t_N=T$ be an uniform partition of $[0,T]$, with $k=t_{i+1}-t_{i}$ for $i = 0,\ldots,N-1$, Set $u^0:= u_0$.
Our semi-discrete time-stepping scheme for \eqref{eq1.1} is defined as follows:

\textbf{Algorithm 1:}

\emph{Step I:} Find $\xi(u^{n-1})\in L^2(\Omega,\cV)$ by solving
\begin{align}\label{eq3.1}
\bigl(\nabla \xi(u^{n-1}),\nabla \phi\bigr) =\bigl(G(u^{n-1}),\nabla \phi\bigr), \ \forall \,\phi \in \cV, \,\, a.s.
\end{align}

\emph{Step II:} Denote $\eta(u^{n-1}):=G(u^{n-1})-\nabla\xi(u^{n-1})$, and find $(u^{n},r^{n})\in\, L^2(\Omega,\cV)\times L^2(\Omega,\cW)$ by solving
\begin{subequations}\label{eq3.2}
\begin{align}\label{eq3.2a}
&\bigl(u^{n},v\bigr) +k\, a\bigl(u^{n}, v\bigr) -k\, d\bigl(v,r^{n}\bigr)+k\, b\bigl(u^{n},u^{n},v\bigr)\\ \nonumber
&\qquad\qquad
=(u^{n-1},v)+ \bigl(\eta(u^{n-1})\Delta W_{n},v\bigr) \quad\forall\, v\in \cV,\,\, a.s.,\\
&d\bigl(u^{n},q\bigr) =0 \qquad\forall\, q\in \cW,\,\, a.s.,  \label{eq3.2b}
\end{align}
\end{subequations}
where $\Delta W_{n}:=W(t_{n})-W(t_{n-1})\thicksim N(0,kQ)$.

\emph{Step III:} Denote $p^{n}:=r^{n}+k^{-1}\xi(u^{n-1})\Delta_{n}W$.

The following lemmas establish some stability results for the
discrete processes $\{(u^n,r^n); 0\leq n\leq N\}$.

\begin{lemma}\label{lem3.3} Let $1\leq p <\infty$ be a natural number.  Assume $u^0\in L^{2^q}(\Omega,\cV_0)$  with $\|u^0\|_{L^2_x}\leq C$.
Then there exists a sequence $\{(u^n,r^n); 1\leq n\leq {N}\}$, which for all $\omega \in \Omega$, solves \textbf{Algorithm 1} and the following stability properties hold:
\begin{align}\label{eq3.3}
\mathbb{E}\Bigl[\max_{1\leq n\leq N}\|u^n\|^{2^p}_{L^2_x}+\nu k\sum^N_{n=1}\|u^n\|^{2^p-1}_{L_x^2}\|\nabla u^n\|^2_{L_x^2}\Bigr] &\leq C.\\\label{eq3.4}
\mathbb{E}\Bigl[\max_{1\leq n\leq N}\|\nabla u^n\|^{2^p}_{L^2_x}+\nu k\sum^N_{n=1}\|\nabla u^n\|^{2^p-1}_{L_x^2}\|\nabla^2 u^n\|^2_{L_x^2}\Bigr] &\leq C.\\\label{eq3.8}
\mathbb{E}\Bigl[\sum^N_{n=1}\|\nabla(u^n-u^{n-1})\|^{2}_{L_x^2}\|\nabla u^n\|^2_{L_x^2}\Bigr] &\leq C,\\\label{eq3.7}
\mathbb{E}\Bigl[\Bigl(\sum^N_{n=1}\|\nabla (u^n-u^n)\|^2_{L^2_x}\Bigr)^4
+\Bigl(k\sum^N_{n=1}\|\nabla^2 u^n\|^{2}_{L_x^2}\Bigr)^4\Bigr]  &\leq C,\\\label{eq3.84}
\mathbb{E}\Bigl[k\sum_{n=1}^{N} \|r^n\|_{L^2_x}^2\Bigr]
&\leq C,\\\label{eq3.5}
\mathbb{E}\Bigl[k\sum_{n=1}^{N} \|\nabla r^n\|_{L^2_x}^2\Bigr]
&\leq C.
\end{align}
\end{lemma}
\begin{proof} Since the proofs of \eqref{eq3.3}--\eqref{eq3.7} were derived in \cite{CP12}. The proofs of \eqref{eq3.84}--\eqref{eq3.5} are similar to \cite{Feng_Qiu18,Feng20}.
For the reader's convenience, we here give their proofs.
 \begin{align}\label{eq3.99}
k\sum_{n=1}^{N}d\bigl(v, r^{n}\bigr)&=\bigl(u^{N},v\bigr)-\bigl(u^{0},v\bigr) +k\sum_{n=1}^{N}\Bigl[a(u^n,v)+b(u^{n},u^n,v)\Bigr] \\ \nonumber
&\quad-\sum_{n=1}^{N} \bigl(\eta(u^{n-1})\Delta W_{n},v\bigr).
\end{align}
Using the Poinc\'{a}re's inequality and the Young's inequality on the right hand side, one finds that
 \begin{align*}
k\sum_{n=1}^{N}d\bigl(v,  r^{n}\bigr)&\leq C\bigl(\|u^{0}\|_{L^2_x} +\|u^n\|_{L^2_x}\bigr)\|\nabla v\|_{L^2_x}
+Ck\sum_{n=1}^{N}\|\nabla u^n\|_{L^2_x}\|\nabla v\|_{L^2_x}\\
&\qquad
+Ck\sum_{n=1}^{N}\bigl(\|u^{n}\|^2_{L^2_x}+\|\nabla u^n\|^2_{L^2_x}\bigr)\|\nabla v\|_{L^2_x}\\
&\qquad\qquad
+C\Bigl\|\sum_{n=1}^{N} \int_{t_{n-1}}^{t_n} \eta(u^{n-1})\, dW(s)\Bigr\|_{L^2_x}\|v\|_{L^2_x}.
\end{align*}
By the inf-sup condition, we get
\begin{align*}
&\beta k\,\sum_{n=1}^{N}\bigl\|  r^n \bigr\|_{L^2_x}
\leq C\bigl(\|u^{0}\|_{L^2_x} +\|u^{N}\|_{L^2_x} \bigr)
+Ck\sum_{n=1}^{N} \|\nabla u^{n}\|_{L^2_x}\\
&\qquad\qquad\qquad\qquad
+Ck\sum_{n=1}^{N}\Bigl(\|u^{n}\|^2_{L^2_x}+\|\nabla u^n\|^2_{L^2_x}\Bigr)
+C\Bigl\|\sum_{n=1}^{N} \int_{t_{n-1}}^{t_n} \eta(u^{n-1})\, dW(s)\Bigr\|_{L^2_x}.
\end{align*}
With the definition of $\eta(u^{n-1})$ and \eqref{eq2.77}, it follows that
\begin{align}\label{eq3.100}
\mathbb{E}\bigl[\|\eta(u^{n-1})\|_{\cK}\bigr]&\leq \mathbb{E}\bigl[\|G(u^{n-1})\|_{L^2_x}+\|\nabla\xi(u^{n-1})\|_{L^2_x}\bigr]\\&\nonumber
\leq C\mathbb{E}\bigr[\|G(u^{n-1})\|_{L^2_x}\bigl].
\end{align}
 Hence, taking the expectation and using It\^{o}'s isometry and \eqref{eq3.3}--\eqref{eq3.4}, which lead to the desired result.

 For \eqref{eq3.5}, setting $v=\nabla r^n$ in \eqref{eq3.2a}, by using the Poinc\'{a}re's inequality and the Young's inequality, one finds that
 \begin{align*}
k\sum_{n=1}^{N}\bigl(\nabla r^{n},\nabla r^n\bigr)&=\bigl(u^{N},\nabla r^n\bigr)-\bigl(u^{0},\nabla r^n\bigr) +k\sum_{n=1}^{N}\Bigl[a(u^n,r^n)+b(u^{n},u^n,\nabla r^n)\Bigr] \\ \nonumber
&\quad-\sum_{n=1}^{N} \bigl(\eta(u^{n-1})\Delta W_{n},\nabla r^n\bigr)\\&
\leq C\bigl(\|u^{0}\|_{L^2_x} +\|u^N\|_{L^2_x}\bigr)\|\nabla r^n\|_{L^2_x}
+C\sum_{n=1}^{N}k\|\nabla^2 u^n\|_{L^2_x}\|\nabla r^n\|_{L^2_x}\\
&\quad
+C \sum_{n=1}^{N}k\bigl(\|\nabla u^{n}\|^2_{L^2_x}+\|\nabla^2 u^{n}\|_{L^2_x}\bigl)\|\nabla r^n\|_{L^2_x}.
\end{align*}
With a standard calculation, we have
 \begin{align*}
 k\sum_{n=1}^{N}\|\nabla r^{n}\|^2_{L^2_x}
&\leq C\bigl(\|u^{0}\|^2_{L^2_x} +\|u^N\|^2_{L^2_x}\bigr)
+Ck\sum_{n=1}^{N}\|\nabla^2 u^n\|^2_{L^2_x} \\
&\quad
+Ck\sum_{n=1}^{N}\bigl(\|\nabla u^{n}\|^2_{L^2_x}+\|\nabla^2 u^{n}\|^2_{L^2_x}\bigl).
\end{align*}
Taking the expectation, using It\^{o}'s isometry, \eqref{eq2.4b},  \eqref{eq3.3} and  \eqref{eq3.4},  we get the desired result. The proof is complete.
\end{proof}


\begin{lemma}\label{lem3.4}
 Assume $u^0\in L^{2^q}(\Omega,H^2)$ for some $1\leq q<\infty$.
Then there exists a sequence $\{u^n; 1\leq n\leq N\}$, which for all $\omega \in \Omega$, solves \textbf{Algorithm 1} and satisfies the following bounds:
\begin{align}\label{eq3.6}
\mathbb{E}\Bigl[\max_{1\leq n\leq N}\|Au^n\|^{2^p}_{L^2_x}
+k\sum^N_{n=1}\|Au^n\|^{2^p-2}_{L_x^2}\|A^{\frac{3}{2}}u^n\|^2_{L_x^2}\Bigr] &\leq C,\\\label{eq3.666}
\mathbb{E}\Bigl[\max_{1\leq n\leq N}\|u^n\|^{2^p}_{-1}
+k\sum^N_{n=1}\|u^n\|^{2^p-2}_{-1}\|u^n\|^2_{L_x^2}\Bigr] &\leq C,\\\label{eq3.6666}
\mathbb{E}\Bigl[\Bigl(\sum^N_{n=1}\|A(u^n-u^n)\|^2_{L^2_x}\Bigr)^4
+\Bigl(k\sum^N_{n=1}\|A^{\frac{3}{2}} u^n\|^{2}_{L_x^2}\Bigr)^4\Bigr]  &\leq C,
\end{align}
where $A:\cV\cap [H^2(D)]^d\to \cV_0$ denotes the Stokes operator (cf. \cite{Temam01}).
\end{lemma}
\begin{proof} For the first assertion, 
taking $v=A^2u^n\in\cV_0$ and $q=0$ in \eqref{eq3.2},  we get
\begin{align}\label{eq3.67}
&\frac{1}{2}\bigl(\|Au^{n}\|^2_{L^2_x}-\|Au^{n-1}\|^2_{L^2_x}+\|A(u^{n}-u^{n-1})\|^2_{L^2_x}\bigr) +k\nu\|A^{\frac{3}{2}}u^n\|^2_{L^2_x}\\ \nonumber
&\qquad\qquad
=kb\bigl(u^{n},u^n,A^2u^n\bigr)+\bigl(A\eta(u^{n-1})W_n,A(u^{n}-u^{n-1})\bigr)
\\ \nonumber
&\quad\qquad\qquad+\bigl(A\eta(u^{n-1})\Delta W_n,A(u^{n-1})\bigr).
\end{align}
 By Lemma 2.1.20 in \cite{K2012} and the Young's inequality,
the first term on the right
hand of \eqref{eq3.67} can be estimated by
\begin{align*}
kb\bigl(u^{n},u^n,A^2u^n\bigr)&\leq Ck\|\nabla^3u^n\|^{\frac{7}{4}}_{L^2_x}\|\nabla u^n\|^{\frac{3}{4}}_{L^2_x}\|u^n\|^{\frac{1}{2}}_{L^2_x}\\&\nonumber
\leq \frac{k\nu}{2}\|\nabla^3u^n\|^{2}_{L^2_x}+Ck\|\nabla u^n\|^{6}_{L^2_x}\|u^n\|^{4}_{L^2_x}\\&\nonumber
\leq \frac{k\nu}{2}\|\nabla^3u^n\|^{2}_{L^2_x}+Ck(\|\nabla u^n\|^{10}_{L^2_x}+\|u^n\|^{10}_{L^2_x}).
\end{align*}
The last term on the right
hand of \eqref{eq3.67} vanishes when taking its expectation. Applying the Young's inequality
 and the tower property for conditional expectations to the second term on the right hand.
 Summing up then leads to
\begin{align}
&\mathbb{E}\Bigl[\|Au^{n}\|^2_{L^2_x}\Bigr]+\mathbb{E}\Bigl[\sum^N_{n=1}\|A(u^{n}-u^{n-1})\|^2_{L^2_x} \Bigr]+\mathbb{E}\Bigl[{\nu}\sum^N_{n=1}k\|A^{\frac{3}{2}}u^n\|^2_{L^2_x}\Bigr]\\ \nonumber
&
\leq C\mathbb{E}\bigl[\|Au^{0}\|^2_{L^2_x}]+C\mathbb{E}\Bigl[\sum^N_{n=1}k\bigl(\|\nabla u^n\|^{10}_{L^2_x}+\|u^n\|^{10}_{L^2_x}\bigr)\Bigr]\\ \nonumber
&\quad
+\mathbb{E}\Bigl[\sum^N_{n=1}\mathbb{E}\bigl[\|P_H\eta(u^{n-1})\|^2_{\mathcal{L}(\cK,H^2)}\|\Delta_nW\|^2_{\cK}|\mathfrak{F}_{t_{n-1}}\bigr]\Bigr]\\ \nonumber
&
\leq C\mathbb{E}\bigl[\|Au^{0}\|^2_{L^2_x}]+C\mathbb{E}\Bigl[\sum^N_{n=1}k\bigl(\|\nabla u^n\|^{10}_{L^2_x}+\|u^n\|^{10}_{L^2_x}\bigr)\Bigr]
+\mathbb{E}\Bigl[\sum^N_{n=1}k\|Au^{n-1}\|^{2}_{L^2_x}\Bigr].
\end{align}
Using Lemma \ref{lem3.3} and the discrete Gronwall's lemma leads to
\begin{align}
\max_{1\leq n\leq N}\mathbb{E}\bigl[\|Au^{n}\|^2_{L^2_x}\bigr] +\mathbb{E}\Bigl[{\nu}\sum^N_{n=1}k\|A^{\frac{3}{2}}u^n\|^2_{L^2_x}\Bigr] \leq C.
\end{align}
To derive the first inequality in \eqref{eq3.6}, using the Young's inequality and Lemma \ref{lem3.3},  one finds that
\begin{align}
\mathbb{E}\bigl[\max_{1\leq n\leq N}\|Au^{n}\|^2_{L^2_x}\bigr] &\leq \mathbb{E}\bigl[\|Au^{0}\|^2_{L^2_x}\bigl]
+C\mathbb{E}\Bigl[\sum^N_{n=1}k\bigl(\|\nabla u^n\|^{10}_{L^2_x}+\|u^n\|^{10}_{L^2_x}\bigr)\Bigr]\\&\nonumber
\quad+C\mathbb{E}\Bigl[\sum^N_{n=1}\|AP_H\eta(u^{n-1})\Delta W_n\|^{2}_{L^2_x}\Bigr]\\&\nonumber\quad
+C\mathbb{E}\Bigl[\max_{1\leq n\leq N}\sum^n_{l=1}(AP_H\eta(u^{l-1})\Delta W_n,Au^{l-1})\Bigr].
\end{align}
The second term on the right hand side may be controlled by Lemma \ref{lem3.3},
 the third term may be estimated by the tower property for conditional expectations,
the fourth term is bounded  with using the
Burkholder-Davis-Gundy inequality. Thus, \eqref{eq3.6} holds for $p=2$. For $p\geq 3$, by the similar line in \cite{BCP12}, we may derive the desired result.
Here we skip it.

For the second assertion, 
taking $v=A^{-1}u^n\in\cV_0$ and $q=0$ in \eqref{eq3.2},  one finds that
\begin{align}\label{eq3.68}
&\frac{1}{2}\bigl(\|u^{n}\|^2_{-1}-\|u^{n-1}\|^2_{-1}+\|u^{n}-u^{n-1}\|^2_{-1}\bigr) +k\nu\|u^n\|^2_{L^2_x}\\ \nonumber
&\qquad\qquad
=kb\bigl(u^{n},u^n,A^{-1}u^n\bigr)+\bigl(A^{-\frac{1}{2}}\eta(u^{n-1})W_n,A^{-\frac{1}{2}}(u^{n}-u^{n-1})\bigr)
\\ \nonumber
&\quad\qquad\qquad+\bigl(A^{-\frac{1}{2}}\eta(u^{n-1})\Delta W_n,A^{-\frac{1}{2}}(u^{n-1})\bigr).
\end{align}
By the Young's inequality and the H\"{o}lder inequality, the first term on the right
hand of \eqref{eq3.68} can be bounded
\begin{align*}
&kb\bigl(u^{n},u^n,A^{-1}u^n\bigr)\\
&\leq Ck\|u^n\|_{L^2_x}\|\nabla u^n\|_{L^2_x}\|u^n-u^{n-1}\|_{-1}+Ck\|u^n\|_{L^2_x}\|\nabla u^n\|_{L^2_x}\|u^{n-1}\|_{-1}\\&
\leq \frac{1}{4}\|u^n-u^{n-1}\|^2_{-1}+Ck^2\|u^n\|^2_{L^2_x}\|\nabla u^n\|^2_{L^2_x}+Ck\|u^n\|_{L^2_x}\|\nabla u^n\|_{L^2_x}\|u^{n-1}\|_{-1}.
\end{align*}
The last term on the right
hand of \eqref{eq3.68} vanishes when taking its expectation.
 Using the Young's inequality
and the tower property for conditional expectations to the second term on the right hand.
 Summing up then leads to
\begin{align}
&\mathbb{E}\Bigl[\|u^{n}\|^2_{-1}\Bigr]+\mathbb{E}\Bigl[\sum^N_{n=1}\|u^{n}-u^{n-1}\|^2_{-1} \Bigr]+\mathbb{E}\Bigl[{\nu}\sum^N_{n=1}k\|u^n\|^2_{L^2_x}\Bigr]\\ \nonumber
&
\leq C\mathbb{E}\Bigl[\sum^N_{n=1}k^2\|u^n\|^2_{L^2_x}\|\nabla u^n\|^2_{L^2_x}\Bigr]+C\mathbb{E}\Bigl[\sum^N_{n=1}k\|u^n\|_{L^2_x}\|\nabla u^n\|_{L^2_x}\|u^{n-1}\|_{-1}\Bigr]\\ \nonumber&
\quad+\mathbb{E}\Bigl[\sum^N_{n=1}\mathbb{E}\bigl[\|P_H\eta(u^{n-1})\|^2_{\mathcal{L}(\cK,H^{-1})}\|\Delta_nW\|^2_{\cK}|\mathfrak{F}_{t_{n-1}}\bigr]\Bigr]\\ \nonumber
&
\leq C\Bigl(\mathbb{E}\Bigl[\sum^N_{n=1}k\|u^n\|^4_{L^2_x}\Bigr]\Bigr)^{\frac{1}{2}}\Bigl(\mathbb{E}\Bigl[\sum^N_{n=1}k\|\nabla u^n\|^4_{L^2_x}\Bigr]\Bigr)^{\frac{1}{2}}\\ \nonumber&
\quad
+C\Bigl(\mathbb{E}\Bigl[\max_{1\leq n\leq N}\|u^n\|^4_{L^2_x}\Bigr]\Bigl)^{\frac{1}{2}}\Bigl(\mathbb{E}\Bigl[\max_{1\leq n\leq N}\|\nabla u^n\|^4_{L^2_x}\Bigr]\Bigl)^{\frac{1}{2}}+\mathbb{E}\Bigl[\sum^N_{n=1}k\|u^{n-1}\|^{2}_{-1}\Bigr].
\end{align}
Using Lemma \ref{lem3.3} and the discrete Gronwall's lemma, one finds that
\begin{align}
\max_{1\leq n\leq N}\mathbb{E}\bigl[\|u^{n}\|^2_{-1}\bigr] +\mathbb{E}\Bigl[{\nu}\sum^N_{n=1}k\|u^n\|^2_{L^2_x}\Bigr] \leq C.
\end{align}
To obtain the first inequality in \eqref{eq3.666}, using the Young's inequality and Lemma \ref{lem3.3},  it follows that
\begin{align}
\mathbb{E}\bigl[\max_{1\leq n\leq N}\|u^{n}\|^2_{-1}\bigr] &\leq \mathbb{E}\bigl[\|u^{0}\|^2_{-1}\bigl]
+C\mathbb{E}\Bigl[\sum^N_{n=1}k\|\nabla u^n\|^{4}_{L^2_x}\Bigr]+C\mathbb{E}\Bigl[\max_{1\leq n\leq N}\|u^n\|^{4}_{L^2_x}\Bigr]\\&\nonumber
\quad+C\mathbb{E}\Bigl[\sum^N_{n=1}\|A^{-\frac{1}{2}}P_H\eta(u^{n-1})\Delta W_n\|^{2}_{L^2_x}\Bigr]\\&\nonumber\quad
+C\mathbb{E}\Bigl[\max_{1\leq n\leq N}\sum^n_{l=1}(A^{-\frac{1}{2}}P_H\eta(u^{l-1})\Delta W_n,A^{-\frac{1}{2}}u^{l-1})\Bigr].
\end{align}
The second term and the third term on the right hand side may be controlled by Lemma \ref{lem3.3},
the fourth term may be estimated by the tower property for conditional expectations,
the fifth term is bounded with using the
Burkholder-Davis-Gundy inequality. Thus, \eqref{eq3.6} holds for $p=2$. For $p\geq 3$, using the similar line in \cite{BCP12}, we skip it.

For the third assertion, using the similar line in \cite{BCP12}, summing over the index $n=1$ in \eqref{eq3.67} and taking the power four, it follows that
\begin{align}
&\Bigl(\sum_{n=1}^N\|A(u^{n}-u^{n-1})\|^2_{L^2_x}\Bigr)^4 +\Bigl(\sum_{n=1}^Nk\nu\|A^{\frac{3}{2}}u^n\|^2_{L^2_x}\Bigr)^4\\ \nonumber
& \leq C\Bigl(\sum^N_{n=1}k\bigl(\|\nabla u^n\|^{10}_{L^2_x}+\|u^n\|^{10}_{L^2_x}\bigr)\Bigr)^4+C\Bigl(\sum^N_{n=1}\|AP_H\eta(u^{n-1})\Delta W_n\|^{2}_{L^2_x}\Bigr)^4
\\ \nonumber
&\quad+C\Bigl(\sum^N_{n=1}(AP_H\eta(u^{n-1})\Delta W_n,Au^{n-1})\Bigr)^4+C\|Au^0\|^8_{L^2_x}.
\end{align}
Taking the expectation, the second term and the third term on the right hand can be bounded as in \cite{BCP12}, it follows that
\begin{align}
&\mathbb{E}\Bigl[\Bigl(\sum_{n=1}^N\|A(u^{n}-u^{n-1})\|^2_{L^2_x}\Bigr)^4 +\Bigl(\sum_{n=1}^Nk\nu\|A^{\frac{3}{2}}u^n\|^2_{L^2_x}\Bigr)^4\Bigr]\\ \nonumber
& \leq C\mathbb{E}\Bigl[\sum^N_{n=1}k\bigl(\|\nabla u^n\|^{40}_{L^2_x}+\|u^n\|^{40}_{L^2_x}\bigr)\Bigr]
+C\mathbb{E}\Bigl[\sum^N_{n=1}k\|Au^{n-1}\|^{8}_{L^2_x}\Bigr]\\& \nonumber\quad
+C\mathbb{E}\bigl[\|Au^0\|^8_{L^2_x}\bigr].
\end{align}
Thanks to Lemma \ref{lem3.3}, the desired result \eqref{eq3.6666} holds. The
proof is complete.
\end{proof}

 Following {\cite{CP12}, for $\epsilon>0$, we define the following sample sets
 \begin{align}\label{eq3.9}
 \Omega_{k}^{\epsilon}=\Bigl\{\omega\in \Omega\bigl|\max_{1\leq n\leq  N}\|\nabla u^n\|^4_{L^2_x}\leq -\epsilon\log\,k\Bigr\}
\end{align}
such that
 \begin{align}\label{eq3.10}
\mathbb{P}(\Omega_{k}^{\epsilon})\geq 1-\dfrac{\mathbb{E}\bigl[\omega\in \Omega\bigl|\max_{1\leq n\leq  N}\|\nabla u^n\|^4_{L^2_x}\bigr]}{-\epsilon\log\,k}\geq 1+\frac{C}{\epsilon\log\,k},
\end{align}
and
 \begin{align}\label{eq3.999}
 \Omega_{\tau}^{\epsilon}=\Bigl\{\omega\in \Omega\bigl|\max_{1\leq n\leq  N}\bigl(\|Au^n\|^2_{L^2_x},\max_{{{t}_{n-1}\leq s\leq{t}_n}}\|u(s)\|^2_{H^2_x}\bigr)\leq -\epsilon\log\,k\Bigr\}
\end{align}
such that
 \begin{align}\label{eq3.1010}
\mathbb{P}(\Omega_{\tau}^{\epsilon})\geq 1-\dfrac{\mathbb{E}\bigl[\omega\in \Omega\bigl|\max_{1\leq n\leq  N}\bigl(\|Au^n\|^2_{L^2_x},\max_{{{t}_{n-1}\leq s\leq{t}_n}}\|u(s)\|^2_{H^2_x}\bigr)\bigr]}{-\epsilon\log\,k}\geq 1+\frac{C}{\epsilon\log\,k}.
\end{align}

By using the similar line in \cite{CP12,Breit21}, the following theorem states and derives the optimal order error estimate for $\{u^n; 1\leq n\leq N\}$  of various spatial norms.

\begin{theorem}\label{thm3.1} Assume that \eqref{eq2.4} holds and that $u_0\in L^{8}(\Omega,\cV_0)$ is an $\mathfrak{F}_0$-measurable random variable.
Let $(u,r)$ be the unique strong solution to \eqref{eq2.7} in the sense of Definition 2, Assume that
\begin{align}\label{eq3.11}
&\mathbb{E}\Bigl[\|u\|_{C^\alpha(0,T;L^4_x)}\Bigr]\leq C,\qquad \mathbb{E}\Bigl[\|u\|_{C^\alpha(0,T;H^{1}_x)}\Bigr]\leq C, \qquad \mathbb{E}\Bigl[\|u\|_{C^\alpha(0,T;H^{2}_x)}\Bigr]\leq C
\end{align}
for some $\alpha\in (0,\frac{1}{2})$. Then, provided that $0<k<k_0$ with $k_0$ sufficiently small, the
following error estimates hold:
\begin{align}\label{eq3.12}
\mathbb{E}\Bigl[\textbf{1}_{\Omega_{k}^{\epsilon}}\Bigl(\max_{1\leq n\leq  N}\|u({t}^n)-{u}^n\|^2_{L^2}
+\nu k\sum_{n=1}^{ N} \|\nabla (u({t}^n)-{u}^n)\|^2_{L^2}\Bigr)\Bigr]&\leq Ck^{2\alpha-\epsilon},\\\label{eq3.112}
\mathbb{E}\Bigl[\textbf{1}_{\Omega_{\tau}^{\epsilon}}\Bigl(\max_{1\leq n\leq  N}\|\nabla(u({t}^n)-{u}^n)\|^2_{L^2}
+\nu k\sum_{n=1}^{ N} \|A(u({t}^n)-{u}^n)\|^2_{L^2}\Bigr)\Bigr]&\leq Ck^{2\alpha-\epsilon},
\end{align}
where $C$ is a positive constant independent of $k$.
\end{theorem}
\begin{proof}
Since the proof of \eqref{eq3.12} was derived in \cite{Breit21}.
Here we prove \eqref{eq3.112}.
For every $n\geq1$, denote $e_u^{n}:={u}(t_n)-{u}^n$,  subtracting \eqref{eq2.7a} from \eqref{eq3.2a} satisfies the following error equation:
\begin{align}\label{eq3.121}
&\bigl(e_u^{n}-e_u^{n-1},v\bigr) +\int^{{t}_n}_{{t}_{n-1}}\,a\bigl(u(s), v\bigr)-a\bigl(u^n, v\bigr)\,ds\\ \nonumber
&\qquad\qquad
=\int^{{t}_n}_{{t}_{n-1}}\,b\bigl(u^{n},u^n,v\bigr)-b\bigl(u(s),u(s),v\bigr)\,ds\\ \nonumber
&\qquad\qquad\quad+\int^{{t}_n}_{{t}_{n-1}}\bigl(\eta(u(s))-\eta(u^{n-1})dW,v\bigr) \,\quad \forall\, v\in \cV_0.
\end{align}
For \eqref{eq3.112}, setting $v=Ae_u^{n}$ in \eqref{eq3.121},  one finds that
\begin{align}\label{eq3.1221}
&\frac{1}{2}\bigl(\|\nabla e_u^{n}\|^2_{L^2_x}-\|\nabla e_u^{n-1}\|^2_{L^2_x}+\|\nabla (e_u^{n}-e_u^{n-1})\|^2_{L^2_x}\bigr)+k\nu\|A e_u^{n}\|^2_{L^2_x}\\ \nonumber&
=\int^{{t}_n}_{{t}_{n-1}}\,a\bigl(u(s)-u({t}^{n}), Ae_u^{n}\bigr)ds
+\int^{{t}_n}_{{t}_{n-1}}\,  b\bigl(u^{n},u^n,e_u^{n}\bigr)-b\bigl(u(s),u(s),Ae_u^{n}\bigr)ds\\ \nonumber&\qquad
+\int^{{t}_n}_{{t}_{n-1}}\bigl(\nabla[\eta(u(s))-\eta(u^{n-1})]dW,\nabla e_u^{n}\bigr)\\ \nonumber&
=D_1+D_2+D_3.
\end{align}
With the Poincar\'{e} inequality and the Young's inequality, the term $D_1$ can be bounded by
 \begin{align}\label{eq3.1231}
D_1&=\int^{{t}_n}_{{t}_{n-1}}\,a\bigl(u(s)-u({t}^{n}), Ae_u^{n}\bigr)ds
\\\nonumber&\leq C\int^{{t}_n}_{{t}_{n-1}}\,\|A(u(s)-u({t}^{n})\|^2_{L^2_x}ds+\frac{k\nu}{16}\|A e_u^{n}\|^2_{L^2_x}
\end{align}
By adding and subtracting suitable terms, we rewrite the nonlinear term $D_2$ as follows:
\begin{align}\label{eq3.1241}
D_2& =\int^{{t}_n}_{{t}_{n-1}}\,  b\bigl(u^{n},u^n,Ae_u^{n}\bigr)-b\bigl(u(s),u(s),Ae_u^{n}\bigr)ds\\&\nonumber
= \int^{{t}_n}_{{t}_{n-1}}\, \Bigl( b\bigl(u^{n}-u(t_{n}),u^n,Ae_u^{n}\bigr)+b\bigl(u(t_{n}-u(s),u^n,Ae_u^{n}\bigr)\\&\nonumber\quad
+b\bigl(u(s),u(t_n)-u(s),Ae_u^{n}\bigr)+b\bigl(u(s),u^{n}-u(t_{n}),Ae_u^{n}\bigr)\Bigr)ds\\&\nonumber
=D_{2,1}+D_{2,2}+D_{2,3}+D_{2,4}.
\end{align}
With the Young's inequality and the Sobolev inequality, we get
\begin{align}\label{eq3.1251}
D_{2,1}&=\int^{{t}_n}_{{t}_{n-1}}\,b\bigl(u^{n}-u(t_{n}),u^n,Ae_u^{n}\bigr)ds\\\nonumber
&\leq \frac{k\nu}{16}\|A e_u^{n}\|^2_{L^2_x}+ Ck\|\nabla e^{n}\|^{2}_{L^2_x}\|A u^n\|^2_{L^2_x},\\\label{eq3.1261}
D_{2,2}&=\int^{{t}_n}_{{t}_{n-1}}\, b\bigl(u(t_{n})-u(s),u^n,Ae_u^{n}\bigr)ds\\&\nonumber
\leq \frac{k\nu}{16}\|Ae_u^{n}\|^2_{L^2_x}+C\int^{{t}_n}_{{t}_{n-1}}\, \|\nabla (u(t_{n})-u(s))\|^2_{L^2_x}\|Au^n\|^2_{L^2_x}ds,\\\label{eq3.1271}
D_{2,3}&=\int^{{t}_n}_{{t}_{n-1}}\, b\bigl(u(s),u(t_n)-u(s),Ae_u^{n}\bigr)ds\\\nonumber&
\leq \frac{k\nu}{16}\|A e_u^{n}\|^2_{L^2_x}+C\int^{{t}_n}_{{t}_{n-1}}\, \|\nabla u(t_{n})-u(s)\|^2_{L^2_x}\|u(s)\|^2_{H^2_x}ds,\\\label{eq3.1272}
D_{2,4}&=\int^{{t}_n}_{{t}_{n-1}}\, b\bigl(u(s),u^{n}-u(t_{n}),Ae_u^{n}\bigr)ds\\\nonumber&
\leq \frac{k\nu}{16}\|A e_u^{n}\|^2_{L^2_x}+Ck\|\nabla e^{n}\|^{2}_{L^2_x}\max_{{{t}_{n-1}\leq s\leq{t}_n}}\|u(s)\|^2_{H^2_x}.
\end{align}
Inserting estimates \eqref{eq3.1231}, \eqref{eq3.1251}--\eqref{eq3.1272} into \eqref{eq3.1221},
applying the summation operator $\sum^N_{n=1}$ and taking the expectation, using \eqref{eq3.11} and Lemma \ref{lem3.3},  it follows that
\begin{align}\label{eq3.1281}
&\mathbb{E}\Bigl[\textbf{1}_{\Omega_{\tau}^{\epsilon}}\Bigl(\max_{1\leq n\leq N}\|\nabla e_u^{n}\|^2_{L^2_x}
+\sum_{n=1}^N\|\nabla (e_u^{n}-e_u^{n-1})\|^2_{L^2_x}+\frac{1}{2}\sum_{n=1}^Nk\nu\|A e_u^{n}\|^2_{L^2_x}\Bigr)\Bigr]\\ \nonumber&
\leq Ck^{2\alpha}+C\mathbb{E}\Bigl[\textbf{1}_{\Omega_{\tau}^{\epsilon}}\Bigl(\sum_{n=1}^Nk\|\nabla e^{n}\|^{2}_{L^2_x}\bigl(\|Au^n\|^2_{L^2_x}
+\max_{{{t}_{n-1}\leq s\leq{t}_n}}\|u(s)\|^2_{H^2_x}\bigr)\Bigr)\Bigr]\\ \nonumber&\qquad
+\mathbb{E}\Bigl[\textbf{1}_{\Omega_{\tau}^{\epsilon}}\Bigl(\sum_{n=1}^N\int^{{t}_n}_{{t}_{n-1}}\bigl(\nabla[\eta(u(s))-\eta(u^{n-1})]dW,\nabla e_u^{n}\bigr)\Bigr)\Bigr].
\end{align}
Using \eqref{eq2.3}, \eqref{eq2.4}, \eqref{eq2.77} and \eqref{eq3.11},  the last term $D_3$ can be estimated by
\begin{align}\label{eq3.1291}
&\mathbb{E}\Bigl[\textbf{1}_{\Omega_{\tau}^{\epsilon}}\Bigl(\sum_{n=1}^N\int^{{t}_n}_{{t}_{n-1}}\bigl(\nabla[\eta(u(s))-\eta(u^{n-1})]dW,\nabla e_u^{n}\bigr)\Bigr]\\&\nonumber
\leq Ck^{2\alpha}+\frac{1}{2}\mathbb{E}\Bigl[\textbf{1}_{\Omega_{\tau}^{\epsilon}}\Bigl(\max_{1\leq n\leq N}\|\nabla e_u^{n}\|^2_{L^2_x}\Bigr)\Bigr]+C\mathbb{E}\Bigl[\textbf{1}_{\Omega_{\tau}^{\epsilon}}\Bigl(\sum_{n=1}^Nk\|\nabla e_u^{n-1}\|^2_{L^2_x}\Bigr)\Bigr]\\&\nonumber\quad
+\frac{1}{2}\mathbb{E}\Bigl[\textbf{1}_{\Omega_{\tau}^{\epsilon}}\Bigl(\sum_{n=1}^N\|\nabla(e_u^{n}-e_u^{n-1})\|^2_{L^2_x}\Bigr)\Bigr].
\end{align}
Then combining \eqref{eq3.1281} with \eqref{eq3.1291} leads to
\begin{align}\label{eq3.1301}
&\mathbb{E}\Bigl[\textbf{1}_{\Omega_{\tau}^{\epsilon}}\Bigl(\max_{1\leq n\leq N}\|\nabla e_u^{n}\|^2_{L^2_x}+\sum_{n=1}^N\|\nabla(e_u^{n}-e_u^{n-1})\|^2_{L^2_x}+\sum_{n=1}^Nk\nu\|Ae_u^{n}\|^2_{L^2_x}\Bigr)\Bigr]\\ \nonumber&
\leq Ck^{2\alpha}+\overline{C}_1\log\,k^{-\epsilon}\mathbb{E}\Bigl[\textbf{1}_{\Omega_{\tau}^{\epsilon}}\Bigl(\sum_{n=1}^Nk\|\nabla e^{n}\|^{2}_{L^2_x}\Bigr)\Bigr]+\overline{C}_2\mathbb{E}\Bigl[\sum_{n=1}^Nk\|\nabla e_u^{n-1}\|^2_{L^2_x}\Bigr]
\\\nonumber&\quad+Ck\underbrace{\mathbb{E}\Bigl[\textbf{1}_{\Omega_{\tau}^{\epsilon}}\Bigl(\sum_{n=1}^N\|\nabla e^{n}\|^{2}_{L^2_x}\|A(u^{n}-u^{n-1})\|^2_{L^2_x}\Bigr)\Bigr]}_{\Theta_1}\\ \nonumber&\quad
+Ck\underbrace{\mathbb{E}\Bigl[\textbf{1}_{\Omega_{\tau}^{\epsilon}}\Bigl(\sum_{n=1}^N\|\nabla e^{n}\|^{2}_{L^2_x}\max_{{{t}_{n-1}\leq s\leq{t}_n}}\|u(s)\|^2_{H^2_x}\Bigr)\Bigr]}_{\Theta_2}.
\end{align}
The terms $\Theta_1$ and $\Theta_2$ may be controlled by Lemmas \ref{lem2.1}, \ref{lem3.3} and \ref{lem3.4}.
If $0< k\leq k_0,\, k_{*}:=\frac{1}{2\overline{C}_1\log\,k_0^{-\epsilon}}<\frac{1}{\overline{C}_1\log\,k_0^{-\epsilon}}$, since $1\leq\frac{1}{1-\overline{C}_1\log(k^{-\epsilon})k}\leq 2$, it follows that
\begin{align}\label{eq3.13001}
&\mathbb{E}\Bigl[\textbf{1}_{\Omega_{k}^{\epsilon}}\Bigl(\max_{1\leq n\leq N}\|\nabla e_u^{n}\|^2_{L^2_x}+\sum_{n=1}^N\|\nabla(e_u^{n}-e_u^{n-1})\|^2_{L^2_x}+\sum_{n=1}^Nk\nu\|Ae_u^{n}\|^2_{L^2_x}\Bigr)\Bigr]\\ \nonumber&
\leq Ck^{2\alpha}+\frac{\overline{C}_1\log\,k^{-\epsilon}}{1-\overline{C}_1k\log\,k^{-\epsilon}}\mathbb{E}\Bigl[\textbf{1}_{\Omega_{k}^{\epsilon}}\Bigl(\sum_{n=1}^Nk\|\nabla e_u^{n-1}\|^{2}_{L^2_x}\Bigr)\Bigr]\\ \nonumber&\quad
+\frac{\overline{C}_2}{1-\overline{C}_1k\log\,k^{-\epsilon}}\mathbb{E}\Bigl[\textbf{1}_{\Omega_{k}^{\epsilon}}\Bigl(\sum_{n=1}^Nk\|e_u^{n-1}\|^2_{L^2_x}\Bigr)\Bigr]
\\ \nonumber&\leq Ck^{2\alpha}+2\bigl({\overline{C}_1\log\,k^{-\epsilon}}+{\overline{C}_2}\bigr)\mathbb{E}\Bigl[\textbf{1}_{\Omega_{k}^{\epsilon}}\Bigl(\sum_{n=1}^Nk\|\nabla e_u^{n-1}\|^{2}_{L^2_x}\Bigr)\Bigr].
\end{align}
By using the discrete Gr{o}nwall inequality, the result \eqref{eq3.112} holds. The proof is complete.
\end{proof}

\begin{remark}
From Lemma \ref{lem2.1}, it is easy to see that the first and second inequalities of the condition \eqref{eq3.11} hold.
Using the similar technique of \cite{Breit21}, the third inequality of the condition \eqref{eq3.11} can be satisfied.
\end{remark}

The last result of this section is stated in the following theorems which give an optimal error estimate for the
pressure $\{r^n; 1\leq n\leq N\}$ and $\{p^n; 1\leq n\leq N\}$.
\begin{theorem}\label{thm3.2}
Let the assumptions of Theorem \ref{thm3.1} be satisfied.
Let $\{r^n; 1\leq n\leq {N}\}$ be the pressure approximation defined by \textbf{Algorithm 1}.
Then the following error estimate holds for $m=1,2,\cdots, N$
\begin{align}\label{eq3.13}
\mathbb{E}\Bigl[\textbf{1}_{\Omega_{k}^{\epsilon}}\Bigl(\Bigl\|\int_0^{t_m}r(s)\,ds-k\sum^{{m}}_{n=1}r^n \Bigr\|^2_{L^2_x}\Bigr)\Bigr]&\leq Ck^{2\alpha-\epsilon},
\end{align}
where $C$ is a positive constant independent of $k$.
\end{theorem}
\begin{proof}
Summing \eqref{eq3.2a} over $1\leq n\leq m(\leq N)$, we get
\begin{align}\label{eq3.14}
&\bigl(u^m,v\bigr) + k\sum_{n=1}^{m} a\bigl(u^n, v\bigr) -k\sum_{n=1}^{m} d\bigl(v,r^n\bigr)+k\sum_{n=1}^{m} b\bigl(u^{n},u^{n},v\bigr)\\ \nonumber
&\qquad
=(u^0,v)+ \sum_{n=1}^{{m}} \bigl(\eta(u^{n-1}) \Delta W_n,v\bigr)
\quad\forall v\in \cV,\,\, a.s.
\end{align}
Subtracting \eqref{eq2.7a} (with $t=t_n$) from \eqref{eq3.14} and noting that $u^0=u(0)$, we obtain
\begin{align}\label{eq3.15}
&d\Bigl(v,k\sum^{{m}}_{n=1}r^n -\int_0^{t_m}r(s)\, ds \Bigr)
=\nu \sum^{{m}}_{n=1} \int^{{t}_n}_{{t}_{n-1}} a\bigl(u(s)-u^n, v)\, ds\\ \nonumber
&\qquad+\sum^{{m}}_{n=1} \int^{{t}_n}_{{t}_{n-1}}  b\bigl(u(s),u(s),v)-b\bigl(u^{n},u^{n},v\bigr)\, ds\\\nonumber
&\qquad +\sum^{{m}}_{n=1} \Bigl(\int^{{t}_n}_{{t}_{n-1}}  \bigl(\eta(u^{n-1})-\eta(u(s))\bigr)\,dW(s),v\Bigr)
+\bigl(u(t^{m})-u^{m},v \bigr) \\ \nonumber
&\quad = \nu\sum^{{m}}_{n=1}\int^{{t}_n}_{{t}_{n-1}} \bigl(\nabla u(s)-\nabla u(t^n)+  \nabla u(t^n)-\nabla u^n,\nabla v\bigr)\, ds\\
\nonumber
&\qquad+\sum^{{m}}_{n=1} \int^{{t}_n}_{{t}_{n-1}} b\bigl(u(s),u(s),v)-b\bigl(u^{n},u^{n},v\bigr)\, ds\\\nonumber
&\qquad +\sum^{{m}}_{n=1} \Bigl(\int^{{t}_n}_{{t}_{n-1}} \bigl(\eta(u^{n-1})-\eta(u(s))\bigr)\,dW(s),v \Bigr)
+\bigl(u(t^{m})-u^{m},v\bigr).
\end{align}
By using the Poincar\'{e} inequality, the H\"{o}lder inequality and the inf-sup condition, one finds that
 \begin{align}\label{eq3.16}
&\beta\Bigl\|k\sum^{{m}}_{n=1}r^n -\int_0^{t_m}r(s)\, ds \Bigr\|_{L^2_x}\\&
\nonumber
\leq C\sum^{{m}}_{n=1}\int^{{t}_n}_{{t}_{n}}\|\nabla u(s)-\nabla u(t^n)\|_{L^2_x}\, ds+  k\sum^{{m}}_{n=1}\|\nabla u(t^n)-\nabla u^n\|_{L^2_x}\\
\nonumber
&\qquad+C\sum^{{m}}_{n=1} \int^{{t}_n}_{{t}_{n-1}} \|u-u(t_{n})\|_{L^4_x}\|\nabla u\|_{L^2_x}+ \|\nabla u(t_{n})\|_{L^2_x}\|\nabla(u-u(t_{n})\|_{L^2_x}\, ds\\\nonumber
&\qquad+C\sum^{{m}}_{n=1}k\Bigl(\|\nabla e^{n}_u\|_{L^2_x}\max_{1\leq l\leq m}\|\nabla u(t_l)\|_{L^2_x}+ \max_{1\leq l\leq m}\|\nabla u^{l}\|_{L^2_x}\|\nabla\,e^{n}_u\|_{L^2_x}\Bigr)\\\nonumber
&\qquad +C\Bigl\|\sum^{{m}}_{n=1}\int^{{t}_n}_{{t}_{n-1}} \bigl(\eta(u^{n-1})-\eta(u(s))\bigr)\,dW(s) \Bigr\|_{L^2_x}
+C\|u(t^{m})-u^{m}\|_{L^2_x}.
\end{align}
Taking the expectation, using \eqref{eq3.11} and Theorem \ref{thm3.1}, it follows that
\begin{align}\label{eq3.17}
\mathbb{E}\Bigl[\textbf{1}_{\Omega_{k}^{\epsilon}}\Bigl(\Bigl\|k\sum^{{m}}_{n=1}r^n &-\int_0^{t_m}r(s)\, ds \Bigr\|_{L^2_x}\Bigr)\Bigr]
\leq Ck^{\alpha}\\\nonumber
&+C\Bigl(\mathbb{E}\Bigl[\textbf{1}_{\Omega_{k}^{\epsilon}}\Bigl(\sum^{{m}}_{n=1}k\|\nabla e^{n}_u\|^2_{L^2_x}\Bigr)\Bigr]\Bigr)^{\frac{1}{2}}\Bigl(\mathbb{E}\max_{1\leq l\leq m}\|\nabla u(t_l)\|^2_{L^2_x}\Bigr)^{\frac{1}{2}}\\\nonumber
&
+C\Bigl(\mathbb{E}\max_{1\leq l\leq m}\|\nabla u^{l}\|^2_{L^2_x}\Bigr)^{\frac{1}{2}}\Bigl(\mathbb{E}\Bigl[\textbf{1}_{\Omega_{k}^{\epsilon}}\Bigl(\sum^{{m}}_{n=1}k\|\nabla\,e^{n}_u\|^2_{L^2_x}\Bigr)\Bigr]\Bigr)^{\frac{1}{2}}
\\&\nonumber+\mathbb{E}\Bigl[\textbf{1}_{\Omega_{k}^{\epsilon}}\Bigl(\Bigl\|\sum^{{m}}_{n=1}\int^{{t}_n}_{{t}_{n-1}} \bigl(\eta(u^{n-1})-\eta(u(s))\bigr)\,dW(s) \Bigr\|_{L^2_x}\Bigr)\Bigr].
\end{align}
By using \eqref{eq2.3}, \eqref{eq2.4}, \eqref{eq2.77} and Theorem \ref{thm3.1},  the last term can be bounded by
\begin{align}\label{eq3.177}
\mathbb{E}\Bigl[\textbf{1}_{\Omega_{k}^{\epsilon}}\Bigl(\Bigl\|\sum^{{m}}_{n=1}\int^{{t}_n}_{{t}_{n-1}} \bigl(\eta(u^{n-1})-\eta(u(s))\bigr)\,dW(s) \Bigr\|_{L^2_x}\Bigr)\Bigr]\leq Ck^{\alpha-\frac{\epsilon}{2}}.
\end{align}
Making use of the Lemmas \ref{lem2.1}, \ref{lem3.3} and Theorem \ref{thm3.1}, the result \eqref{eq3.13} holds.  The proof is complete.
\end{proof}

\begin{theorem}\label{thm3.3}
Let the assumptions of Theorem \ref{thm3.1} be satisfied.
Let $\{p^n; 1\leq n\leq {N}\}$ be the pressure approximation defined by \textbf{Algorithm 1}.
Then the following error estimate holds for $m=1,2,\cdots, N$
\begin{align}\label{eq3.18}
\mathbb{E}\Bigl[\textbf{1}_{\Omega_{k}^{\epsilon}}\Bigl(\Bigl\|\int_0^{t_m}p(s)\,ds-k\sum^{m}_{n=1}p^n \Bigr\|^2_{L^2_x} \Bigr)\Bigr]&\leq Ck^{2\alpha-\epsilon},
\end{align}
where $C$ is a positive constant independent of $k$.
\end{theorem}

%
\section{Fully discrete mixed finite element scheme}\label{sec-4}
In this section we propose and analyze a fully discrete time-stepping scheme for the mixed formulation
\eqref{eq2.7}. The error estimates in strong norms
for both the velocity and pressure approximations are obtained.  Furthermore,
we derive strong optimal $H^1$ convergence first, and then obtain $L^2$ convergence of the fully discrete scheme
with the negative norm technique.

 Suppose that $\cT_h$ is a quasi-uniform family of
triangulation of the periodic domain $D\subset \mathbb{R}^2$. We define three
finite element spaces as follows:
\begin{align*}
\cV_h &=\bigl\{v_h\in H^{1}_{per}(D)^2;\, v_h|_K\in P_{2}(K)^2\,\,\forall\, K\in  \cT_h\bigr\},\\
\cW_h &=\bigl\{q_h\in L^2_{per}(D);\, q_h|_K\in P_1(K)\,\,\forall\, K\in \cT_h\bigr\},\\
\cS_h &=\bigl\{w_h\in H^1_{per}(D)^2;\, w_h|_K\in P_{1}(K)^2\,\,\forall\, K\in  \cT_h\bigr\},
\end{align*}
where $P_l(K)$ ($l=1,2$) denotes the set of polynomials of degree less than or equal to $l$
over the element $K\in \cT_h$.

In addition, we consider the weakly discrete divergence-free subspace $\cV_{0h}\subset\cV_{h}$
\begin{equation}\label{weakly_divergence-free}
\cV_{0h}=\Bigl\{v_h\in\cV_h;\, d(q_h,v_h)=0, \, \forall\, q_h \in\ \cW_h \Bigr\}.
\end{equation}

As it is noted \cite{Brezzi_Fortin91} that the finite element space pair $(\cV_h,\cW_h)$
is stable in the sense that the following discrete inf-sup condition holds, i.e., there
exists an $h$-independent positive constant $\gamma$ such that
\begin{equation}\label{inf-sup-TH}
\sup_{v_h\in \cV_h} \frac{ d(v_h,q_h)}{\|\nabla v_h\|_{L^2_x}} \geq \gamma \|q_h\|_{L^2_x}
\qquad \forall\, q_h\in \cW_h.
\end{equation}

We define the $L^2(D)$ projections $\rho_h: L_{per}^2(D)\rightarrow \cW_h$, $\Pi_h: L_{per}^2(D)^2\rightarrow \cV_h$ and
the $L^2$ Ritz-projection $\sigma_h:H^1_{per}(D)^2\rightarrow S_h$ such that
\begin{align*}
 (\varphi-\rho_h\varphi,\psi_h)&=0, \qquad \forall\,\psi_h\in \cW_h,\\
 (v-\Pi_hv,w_h)&=0, \qquad \forall\, w_h\in \cV_h,\\
 (\nabla(\phi-\sigma_h\phi),\nabla\chi_h)&=0, \qquad \forall\, \chi_h\in S_h.
\end{align*}
The following approximation properties are well-known \cite{Girault_Raviart86,Brezzi_Fortin91,Ern_Guermond04,Falk08}
\begin{align}\label{eq4.23}
\|\varphi-\rho_h\varphi\|_{L^2_x} +h\|\nabla(\varphi-\rho_h\varphi)\|_{L^2_x} &\leq Ch^{s}\,\|\varphi\|_{H_x^s} \,\,\forall\, \varphi\in H_{per}^s(D),\\\label{eq4.24}
\|v-\Pi_h v\|_{L_x^2} +h\|\nabla(v-\Pi_h v)\|_{L_x^2} &\leq Ch^{s}\,\|v\|_{H_x^s} \,\,\forall\, v\in H_{per}^s(D)^2,\\\label{eq4.25}
\|\phi-\sigma_h\phi\|_{L^2_x}+ h\|\nabla(\phi-\sigma_h\phi)\|_{L^2_x}&\leq Ch^s\,\|\phi\|_{H_x^s} \,\, \forall\, \phi\in H_{per}^s(D)^2/{\mathbb{R}},
\end{align}
 where $C$ is a positive constant independent of $h$.

Our fully discrete finite element algorithm for \eqref{eq2.7} is defined as follows.

\textbf{Algorithm 2:}

Set $u^0_h\in L^2(\Omega,\cV_h)$,
for $n=1,\ldots,N$, we define the following steps:

\emph{Step I:} Find $\xi(u^{n-1}_h)\in L^2(\Omega,\cS_h)$ by solving
\begin{align}\label{eq4.26}
\bigl(\nabla \xi(u_h^{n-1}),\nabla \phi\bigr) =\bigl(G(u_h^{n-1}),\nabla \phi\bigr), \ \forall \, \phi\,\in \cS_h \,\, a.s.
\end{align}

\emph{Step II:} Denote $\eta(u_{h}^{n-1}):=G(u_{h}^{n-1})-\nabla\xi(u_{h}^{n-1})$, and find $(u_{h}^{n},r_{h}^{n})\in L^2(\Omega,\cV_h)\times L^2(\Omega,\cW_h)$ by solving
\begin{subequations}\label{eq4.27}
\begin{align}\label{eq4.27a}
&\bigl(u^{n}_h,v\bigr) +k\,a\bigl(u_h^{n}, v\bigr) -k\, d\bigl(v,r_h^{n}\bigr)+k\, b\bigl(u_h^{n},u_h^{n},v\bigr)\\ \nonumber
&\qquad\qquad
=(u^{{n-1}}_h,v)+ \bigl(\eta(u_h^{n-1})\Delta W_{n},v\bigr) \quad\forall\, v\in\, \cV_h,\,\, a.s.,\\\label{eq4.27b}
&d\bigl(u^{n}_h,q\bigr) =0 \qquad\forall\, q\,\in\, \cW_h,\,\, a.s.
\end{align}
\end{subequations}

\emph{Step III:} Denote $p_h^{n}:=r_h^{n}+k^{-1}\xi(u_h^{n-1})\Delta_{n}W$.

We now give the following stabilities for $u_h^{n}$ and $r_h^{n}$, but omit their
proofs because they are similar to semi-discrete scheme given in \cite{CP12,Feng_Qiu18}.

\begin{lemma}\label{lem4.5}
Let $1\leq q <\infty$ be a natural number. Assume $u^0_h\in L^{2^q}(\Omega,\cV_{0h})$  with $\|u^0_h\|_{L^2_x}\leq C$.
Let $\{ (u^n_h, r_h^n,p_h^n); 1\leq n\leq N \}$ be a solution to \textbf{Algorithm 2}, then there hold
\begin{subequations}\label{eq4.28}
\begin{align}\label{eq4.28a}
\mathbb{E}\Bigl[\max_{1\leq n\leq N}\|u^n_h\|^{2^p}_{L^2_x}+\nu k\sum^N_{n=1}\|u^n_h\|^{2^p-1}_{L_x^2}\|\nabla u^n_h\|^2_{L_x^2}\Bigr] &\leq C,\\\label{eq3.28b}
\mathbb{E}\Bigl[ k\sum_{n=1}^N \|r_h^n \|_{L^2_x}^2\Bigr]&\leq C.
\end{align}
\end{subequations}
\end{lemma}
\begin{lemma}\label{lem4.6}
Let $1\leq p <\infty$ be a natural number. Assume $u^0_h\in L^{2^q}(\Omega,\cV_{0h})$  with $\|u^0_h\|_{L^2_x}\leq C$.
Then there exists a sequence $\{u^n_h\}_{n\geq1}$ of $\cV$-valued random variables, which for all $\omega \in \Omega$, solves\textbf{ Algorithm 2} and has the following stability estimates:
\begin{subequations}\label{eq4.29}
\begin{align}\label{eq4.29a}
\mathbb{E}\Bigl[\max_{1\leq n\leq N}\|\nabla u_h^n\|^{2^p}_{L^2_x}+\nu k\sum^N_{n=1}\|\nabla u_h^n\|^{2^p-1}_{L_x^2}\|\nabla^2u_h^n\|^2_{L_x^2}\Bigr] &\leq C,\\\label{eq4.29b}
\mathbb{E}\Bigl[\sum^N_{n=1}\|\nabla(u^n_h-u_h^{n-1})\|^{2}_{L_x^2}\|\nabla u_h^n\|^2_{L_x^2}\Bigr] &\leq C,\\\label{eq4.29c}
\mathbb{E}\Bigl[\Bigl(\sum^N_{n=1}\|\nabla (u_h^n-u_h^{n-1})\|^2_{L^2_x}\Bigr)^4
+\Bigl(k\sum^N_{n=1}\|\nabla^2 u_h^n\|^{2}_{L_x^2}\Bigr)^4\Bigr]  &\leq C.
\end{align}
\end{subequations}
\end{lemma}

For $\epsilon>0$, we introduce the sample set
 \begin{align}\label{eq4.99}
 \Omega_{h}^{\epsilon}=\Bigl\{\omega\in \Omega\bigl|\max_{1\leq n\leq  N}\|\nabla u^n\|^4_{L^2_x}+\|u_h^n\|^2_{L^2_x}\leq -\epsilon\log(h^2+k)\Bigr\}
\end{align}
such that
 \begin{align}
\mathbb{P}(\Omega_{h}^{\epsilon})\geq 1-\dfrac{\mathbb{E}\bigl[\omega\in \Omega\bigl|\max_{1\leq n\leq  N}\bigl(\|\nabla u^n\|^4_{L^2_x}+\|u_h^n\|^2_{L^2_x}\bigr)\bigr]}{-\epsilon\log(h^2+k)}\geq 1+\frac{C}{\epsilon\log(h^2+k)}.
\end{align}

We are now in a position to state and prove the first main theorem of this section.
\begin{theorem}\label{thm4.4}
Set $u^0=u_0$ and let  $\{u^n; 1\leq n\leq N\}$ and
$\{u^n_h; 1\leq n\leq N\}$ be the solutions of \textbf{ Algorithm 1} and \textbf{ Algorithm 2},
respectively.
Then, provided that $0<k<k_0$ and $0<h<h_0$ with $k_0$ and $h_0$ sufficiently small, the following error estimate holds:
\begin{align}\label{eq4.30}
& \mathbb{E}\bigl[ \textbf{1}_{\Omega_{h}^{\epsilon}}\Bigl(\max_{1\leq n\leq N}\|{u}^n-{u}^n_{h}\|^2_{L^2} +k\sum_{n=1}^{N}\|\nabla ({u}^n-{u}^n_{h})\|^2_{L^2}\Bigr)\Bigr] \leq C(h^{2-2\epsilon}+k^{1-\epsilon}).
\end{align}
\end{theorem}
\begin{remark} The proof of Theorem \ref{thm4.4} is similar line to the proofs of \cite{CP12,Breit21}. But we use different indicator function \eqref{eq4.99}.
\end{remark}
\begin{proof}
for every $n\geq1$, let $e_u^{n,h}:={u}^n-{u}^n_{h}$ and $e_r^{n,h}:={r}^n-{r}^n_{h}$, it is easy to
check that $(e_u^{n,h},e_r^{n,h})$ satisfies the following error equations:
\begin{subequations}\label{eq4.31}
\begin{align}\label{eq4.31a}
\bigl(e_u^{n,h}-e_u^{n-1,h},v_h\bigr) &+k\,a\bigl(e_u^{n,h},v_h\bigr) +k\,d\bigl(v_h,e_r^{n,h}\bigr) \\ \nonumber
&
+k\,b\bigl({u}^{n},{u}^{n},v_h\bigr)-k\,b\bigl({u}_{h}^{n},{u}_{h}^{n},v_h\bigr)\\ \nonumber
&=\bigl([{\eta}(u^{n-1})-{\eta}(u^{n-1}_h)]\Delta W_{n},v_h \bigr),\quad \forall v_h\,\in\,\cV_h,\\\label{eq4.31b}
d\bigl(e_u^{n,h},q_h\bigr)&=0,\quad \quad \forall q_h\,\in\,\cW_h.
\end{align}
\end{subequations}
Setting $v_h=\Pi_he_u^{n,h}$ and $q_h=\rho_he_r^{n,h}$ in \eqref{eq4.31}, we have
\begin{align}\label{eq4.32}
\bigl(e_u^{n,h}-e_u^{n-1,h},\Pi_he_u^{n,h}\bigr) &+k\,a\bigl(e_u^{n,h},\Pi_he_u^{n,h}\bigr) -k\,d\bigl(\Pi_he_u^{n,h},e_r^{n,h}\bigr) \\ \nonumber
&
+k\,b\bigl({u}^{n},{u}^{n},\Pi_he_u^{n,h}\bigr)-k\,b\bigl({u}_{h}^{n},{u}_{h}^{n},\Pi_he_u^{n,h}\bigr)\\ \nonumber
&=\bigl([{\eta}(u^{n-1})-{\eta}(u^{n-1}_h)]\Delta W_{n},\Pi_he_u^{n,h} \bigr)
\end{align}
By using the identity $a\cdot(a-b)=\frac{1}{2}(|a|^2-|b|^2+|a|^2-|b|^2)$,
we gain
\begin{align}\label{eq4.33}
&\frac{1}{2}\bigl(\|\Pi_he_u^{n,h}\|^2_{L^2_x}-\|\Pi_he_u^{n-1,h}\|^2_{L^2_x}+\|\Pi_he_u^{n,h}-\Pi_he_u^{n-1,h}\|^2_{L^2_x}\bigr)\\ \nonumber
&\quad +k\,\nu\|\nabla \Pi_he_u^{n,h}\|^2_{L^2_x}
=k\,a\bigl(u^n-u^{n}_{h},\Pi_he_u^{n,h}\bigr)+k\,d\bigl(\Pi_he_u^{n,h},e_p^{n,h}\bigr) \\ \nonumber
&\quad+k\,b\bigl({u}_{h}^{n},{u}_{h}^{n},\Pi_he_u^{n,h}\bigr)-k\,b\bigl({u}^{n},{u}^{n},\Pi_he_u^{n,h}\bigr)\\ \nonumber
&\quad+\bigl([{\eta}(u^{n-1})-{\eta}(u^{n-1}_h)]\Delta W_{n},\Pi_he_u^{n,h} \bigr)\\ \nonumber&
=\sum^4_{i=1}I_i.
\end{align}
For terms $I_1$ and $I_2$, thanks to the Young's inequality, \eqref{eq4.23} and \eqref{eq4.24}, we obtain
\begin{align}\label{eq4.34}
I_1 &
\leq \frac{\nu k}{8}\|\nabla \Pi_he_u^{n,h}\|^2_{L^2_x}+C{kh^2}\|\nabla^2{u}^{n}\|^2_{L^2_x},\\ \label{eq4.35}
I_2  
&  \leq \frac{\nu k}{8}\|\nabla \Pi_he_u^{n,h}\|^2_{L^2_x}+Ckh^2\|\nabla r^{n}\|_{L^2_x}^2.
\end{align}
For nonlinear term $I_3$, we rewrite as follows:
\begin{align*}
I_3&=-kb\bigl({u}^{n}-\Pi_h{u}^{n},u^{n},\Pi_he_u^{n,h}\bigr)-kb\bigl(\Pi_he_u^{n,h},\Pi_he_u^{n,h},{u}^{n}\bigr)\\&\quad
-kb\bigl(\Pi_he_u^{n,h},\Pi_he_u^{n,h},{u}^{n}-\Pi_h{u}^{n}\bigr)+kb\bigl(\Pi_hu^{n},\Pi_he_u^{n,h},{u}^{n}-\Pi_h{u}^{n}\bigr)\\&
=\sum^4_{i=1}I_{3,i}.
\end{align*}
Using the Poincar\'{e} inequality, the Young's inequality, the embedding inequality and \eqref{eq4.24}, one finds that
\begin{align}\label{eq4.36}
I_{3,1}&\leq Ck\|{u}^{n}-\Pi_h{u}^{n}\|^{\frac{1}{2}}_{L^2_x}\|\nabla({u}^{n}-\Pi_h{u}^{n})\|^{\frac{1}{2}}_{L^2_x}\|\nabla u^{n}\|_{L^2_x}\|\nabla\Pi_he_u^{n,h}\|_{L^2_x}\\&\nonumber
\leq \frac{k\nu}{16}\|\nabla \Pi_he_u^{n,h}\|^2_{L^2_x}+Ckh^3\|\nabla^2{u}^{n}\|^2_{L^2_x}\|\nabla{u}^{n}\|^2_{L^2_x},\\\label{eq4.37}
I_{3,2}&\leq k\|\Pi_he_u^{n,h}\|^{\frac{1}{2}}_{L^2_x}\|\nabla \Pi_he_u^{n,h}\|^{\frac{1}{2}}_{L^2_x}\|\nabla\Pi_he_u^{n,h}\|_{L^2_x}\|u^{n}\|_{L^4_x}\\&\nonumber
\leq \frac{k\nu}{16}\|\nabla \Pi_he_u^{n,h}\|^2_{L^2_x}+Ck\|\Pi_he_u^{n,h}\|^2_{L^2_x}\|\nabla{u}^{n}\|^4_{L^2_x},\\\label{eq4.38}
I_{3,3}&\leq k\|\Pi_he_u^{n,h}\|^{\frac{1}{2}}_{L^2_x}\|\nabla \Pi_he_u^{n,h}\|^{\frac{1}{2}}_{L^2_x}
\|\nabla \Pi_he_u^{n,h}\|_{L^2_x}\|{u}^{n}-\Pi_h{u}^{n}\|_{L^4_x}\\&\nonumber
\leq \frac{k\nu}{16}\|\nabla \Pi_he_u^{n,h}\|^2_{L^2_x}+Ckh^2\|e_u^{n,h}\|^{2}_{L^2_x}\|\nabla{u}^{n}\|^{4}_{L^2_x},\\\label{eq4.39}
I_{3,4}&\leq k\|\nabla \Pi_hu^{n}\|_{L^2_x}\|\nabla \Pi_he_u^{n,h}\|_{L^2_x}\|{u}^{n}-\Pi_h{u}^{n}\|_{L^4_x}\\&\nonumber
\leq \frac{k\nu}{16}\|\nabla \Pi_he_u^{n,h}\|^2_{L^2_x}+Ckh^3\|\nabla u^{n}\|^{2}_{L^2_x}\|\nabla^2{u}^{n}\|^2_{L^2_x}.
\end{align}
Inserting estimates \eqref{eq4.34}--\eqref{eq4.39} into \eqref{eq4.33}, we arrive at
\begin{align}\label{eq4.40}
&\frac{1}{2}\bigl(\|\Pi_he_u^{n,h}\|^2_{L^2_x}-\|\Pi_he_u^{n-1,h}\|^2_{L^2_x}+\|\Pi_he_u^{n,h}-\Pi_he_u^{n-1,h}\|^2_{L^2_x}\bigr)\\ \nonumber&\quad
+\frac{k\nu}{2}\|\nabla\Pi_he_u^{n,h}\|^2_{L^2_x}
\leq C{kh^2}\|\nabla^2{u}^{n}\|^2_{L^2_x}+Ckh^2\|\nabla r^{n}\|_{L^2_x}^2\\ \nonumber&\qquad
+Ckh^3\|\nabla^2{u}^{n}\|^2_{L^2_x}\|\nabla{u}^{n}\|^2_{L^2_x}
 +Ck\|\Pi_he_u^{n,h}\|^2_{L^2_x}\|\nabla{u}^{n}\|^4_{L^2_x}\\ \nonumber&\qquad\quad
 +Ckh^2\|e_u^{n,h}\|^{2}_{L^2_x}\|\nabla{u}^{n}\|^{2}_{L^2_x}
+Ckh^3\|\nabla u^{n}\|^{2}_{L^2_x}\|\nabla^2{u}^{n}\|^2_{L^2_x}\\ \nonumber &\qquad\qquad
+\bigl([{\eta}(u^{n-1})-{\eta}(u^{n-1}_h)]\Delta W_{n},\Pi_he_u^{n,h} \bigr)
\end{align}
Taking the expectation and applying the summation operator $\sum^N_{n=1}$, one finds that
\begin{align}\label{eq4.41}
&\mathbb{E}\Big[\textbf{1}_{\Omega_{h}^{\epsilon}}\Bigl(\frac{1}{2}\|\Pi_he_u^{n,h}\|^2_{L^2_x}+\frac{1}{2}\sum^N_{n=1}\|\Pi_he_u^{n,h}-\Pi_he_u^{n-1,h}\|^2_{L^2_x}+\sum^N_{n=1}\frac{k\nu}{2}\|\nabla \Pi_he_u^{n,h}\|^2_{L^2_x}\Bigr)\Bigr]\\ \nonumber&
\leq \frac{1}{2}\mathbb{E}\big[\|\Pi_he_u^{0,h}\|^2_{L^2_x}\bigr]+C{h^2}\mathbb{E}\Big[k\sum^N_{n=1}\|\nabla^2{u}^{n}\|^2_{L^2_x}\Bigr]+Ch^2\mathbb{E}\Big[k\sum^N_{n=1}\|\nabla r^{n}\|_{L^2_x}^2\Bigr]\\ \nonumber&\quad
+Ch^3\underbrace{\mathbb{E}\Big[k\sum^N_{n=1}\|\nabla^2{u}^{n}\|^2_{L^2_x}\|\nabla{u}^{n}\|^2_{L^2_x}\Bigr]}_{\Lambda_1}
+Ck\underbrace{\mathbb{E}\Big[\sum^N_{n=1}\|e_u^{n,h}\|^2_{L^2_x}\|\nabla({u}^{n}-{u}^{n-1})\|^4_{L^2_x}\Bigr]}_{\Lambda_2}\\ \nonumber&\quad
+Ch^2\underbrace{\mathbb{E}\Big[k\sum^N_{n=1}\|e_u^{n,h}\|^{2}_{L^2_x}\|\nabla{u}^{n}\|^{4}_{L^2_x}\Bigr]}_{\Lambda_3}
+Ch^3\underbrace{\mathbb{E}\Big[k\sum^N_{n=1}\|\nabla u^{n}\|^{2}_{L^2_x}\|\nabla^2{u}^{n}\|^2_{L^2_x}\Bigr]}_{\Lambda_4}\\ \nonumber&\quad
+C\log(h^2+k)^{-\epsilon}\mathbb{E}\Big[\textbf{1}_{\Omega_{h}^{\epsilon}}\Bigl(\sum^N_{n=1}k\|\Pi_he_u^{n,h}\|^2_{L^2_x}\Bigr)\Bigr]\\ \nonumber&\quad
+\mathbb{E}\Big[\textbf{1}_{\Omega_{h}^{\epsilon}}\Bigl(\sum^N_{n=1}\bigl([{\eta}(u^{n-1})-{\eta}(u^{n-1}_h)]\Delta W_{n},\Pi_he_u^{n,h} \bigr)\Bigr)\Bigr].
\end{align}
Now we explain how to estimate in expectation for $\Lambda_i$ $(i=1,\ldots,4)$. Making use of the Lemmas \ref{lem3.3} and \ref{lem4.5}, the terms  $\Lambda_i$ $(i=1,3)$ are uniformly bounded
\begin{align*}
{\Lambda_1}&
\leq  \Bigl(\mathbb{E}\max_{1\leq m\leq N}\|\nabla{u}^{m}\|^4_{L^2_x}\Bigr)^{\frac{1}{2}}\Bigl(\mathbb{E}\sum^N_{n=1}k\|\nabla^2{u}^{n}\|^4_{L^2_x}\Bigr)^{\frac{1}{2}},\\
{\Lambda_3}&
\leq  \Bigl(\mathbb{E}\max_{1\leq m\leq N}\|u^{m}\|^4_{L^2_x}\Bigr)^{\frac{1}{2}}\Bigl(\mathbb{E}\sum^N_{n=1}k\|\nabla{u}^{n}\|^8_{L^2_x}\Bigr)^{\frac{1}{2}}\\&\nonumber\quad
+ \Bigl(\mathbb{E}\max_{1\leq m\leq N}\|u^{m}_h\|^4_{L^2_x}\Bigr)^{\frac{1}{2}}\Bigl(\mathbb{E}\sum^N_{n=1}k\|\nabla{u}^{n}\|^8_{L^2_x}\Bigr)^{\frac{1}{2}}.
\end{align*}
About the term  $\Lambda_2$, using the Lemmas \ref{lem3.3} and \ref{lem4.5}, we have
\begin{align*}
&\mathbb{E}\Big[\sum^N_{n=1}\|e_u^{n,h}\|^2_{L^2_x}\|\nabla({u}^{n}-{u}^{n-1})\|^4_{L^2_x}\Bigr]\\
&
\leq  \Big(\mathbb{E}\max_{1\leq m\leq N}\|u^m\|^4_{L^2_x}\Bigr)^{1/2}\Big(\mathbb{E}\Big[\sum^N_{n=1}\|\nabla({u}^{n}-{u}^{n-1})\|^8_{L^2_x}\Bigr]\Bigr)^{1/2}\\&\quad+
\Big(\mathbb{E}\max_{1\leq m\leq N}\|u^m_h\|^4_{L^2_x}\Bigr)^{1/2}\Big(\mathbb{E}\Big[\sum^N_{n=1}\|\nabla({u}^{n}-{u}^{n-1})\|^8_{L^2_x}\Bigr]\Bigr)^{1/2},
\end{align*}
and the term  $\Lambda_4$ is uniformly bounded as follows:
\begin{align*}
{\Lambda_4}&
\leq \Bigl(\mathbb{E}\max_{1\leq m\leq N}\|\nabla u^{m}\|^4_{L^2_x}\Bigr)^{\frac{1}{2}}\Bigl(\mathbb{E}\sum^N_{n=1}k\|\nabla^2{u}^{n}\|^4_{L^2_x}\Bigr)^{\frac{1}{2}}.
\end{align*}
For term $I_4$, using It\^{o}'s isometry and the Young's inequality, we have
\begin{align}\label{eq4.42}
&\mathbb{E}\Bigl[\textbf{1}_{\Omega_{h}^{\epsilon}}\Bigl(v\sum_{n=1}^{N} \bigl([{\eta}(u^{n-1})-{\eta}(u^{n-1}_h)]\Delta W_{n},\Pi_he_u^{n,h} \bigr) \Bigr)\Bigr]\\
\nonumber
&
=\mathbb{E}\Bigl[\textbf{1}_{\Omega_{h}^{\epsilon}}\Bigl(\sum_{n=1}^{N} \bigl([{\eta}(u^{n-1})-{\eta}(u^{n-1}_h)]\Delta W_{n},\Pi_he_u^{n,h}-\Pi_he_u^{n-1,h} \bigr) \Bigr)\Bigr]\\ \nonumber&
\leq \mathbb{E}\Bigl[\textbf{1}_{\Omega_{h}^{\epsilon}}\Bigl(k\sum_{n=1}^{N} \bigl\|[{\eta}(u^{n-1})-{\eta}(u^{n-1}_h)]\Delta W_{n}\|_{L^2_x}\|\Pi_he_u^{n,h}-\Pi_he_u^{n-1,h}\|_{L^2_x} \Bigr)\Bigr]\\ \nonumber&
\leq \frac14 \mathbb{E}\Bigl[\textbf{1}_{\Omega_{h}^{\epsilon}}\Bigl(\|\Pi_he_u^{n,h}-\Pi_he_u^{n-1,h}\|^2_{L^2_x}\Bigr)\Bigr]+ \mathbb{E}\Bigl[\textbf{1}_{\Omega_{h}^{\epsilon}}\Bigl(k\sum_{n=1}^{N} \bigl\|{\eta}(u^{n-1})-{\eta}(u^{n-1}_h)\|^2_{L^2_x}\Bigr) \Bigr].
\end{align}
With the definition of ${\eta}$, and using \eqref{eq2.3}, \eqref{eq2.4}, \eqref{eq2.77}, \eqref{eq3.11} and  \eqref{eq4.24}, one finds that
\begin{align}\label{eq4.46}
&\quad\|{\eta}(u^{n-1})-{\eta}(u_h^{n-1})\|^2_{L^2_x}\\&\nonumber \leq C\|e_u^{n-1}\|_{L^2_x}^2+Ch^2\|{\eta}(u^{n-1})\|_{H^2_x}^2\\&\nonumber
\leq Ch^2\|\nabla\cdot G({u}^{n-1})\|_{L^2}^2+C\|e_u^{n-1}\|_{L^2_x}^2\\&\nonumber
\leq Ch^2\|\nabla {u}^{n-1}\|_{L^2_x}^2+Ch^4\|\nabla^2 {u}^{n-1}\|_{L^2_x}^2+C\|\Pi_he_u^{n-1,h}\|_{L^2_x}^2.
\end{align}
Combining \eqref{eq4.42}--\eqref{eq4.46} into \eqref{eq4.41}, we get
\begin{align}\label{eq4.48}
&\mathbb{E}\Big[\textbf{1}_{\Omega_{h}^{\epsilon}}\Bigl(\frac{1}{2}\|\Pi_he_u^{n,h}\|^2_{L^2_x}+\frac{1}{4}\sum^N_{n=1}\|\Pi_he_u^{n,h}-\Pi_he_u^{n-1,h}\|^2_{L^2_x}+\sum^N_{n=1}\frac{k\nu}{2}\|\nabla e_u^{n,h}\|^2_{L^2_x}\Bigr)\Bigr]\\ \nonumber&
\leq C(h^2+k)+C_3\log(h^2+k)^{-\epsilon}\mathbb{E}\big[\sum^N_{n=1}k\|\Pi_he_u^{n,h}\|^2_{L^2_x}\bigr]
\\ \nonumber&\quad+C_4\mathbb{E}\Bigl[\sum_{n=1}^Nk\|\Pi_he_u^{n-1,h}\|_{L^2_x}^2\Bigr].
\end{align}
If $0<h<h_0$ and $0< k\leq k_0, \,k^{*}:=\frac{1}{2C_3\log(h_0^2+k_0)^{-\epsilon}}<\frac{1}{C_3\log(h_0^2+k_0)^{-\epsilon}}$, since $1\leq\frac{1}{1-C_3k\log(h^2+k)^{-\epsilon}}\leq 2$, it follows that
\begin{align}
&\mathbb{E}\Big[\textbf{1}_{\Omega_{h}^{\epsilon}}\Bigl(\frac{1}{2}\|\Pi_he_u^{n,h}\|^2_{L^2_x}+\frac{1}{4}\sum^N_{n=1}\|\Pi_he_u^{n,h}-\Pi_he_u^{n-1,h}\|^2_{L^2_x}+\sum^N_{n=1}\frac{k\nu}{2}\|\nabla e_u^{n,h}\|^2_{L^2_x}\Bigr)\Bigr]\\ \nonumber&
\leq C(h^2+k)+\frac{C_3\log{(h^2+k)^{-\epsilon}}}{1-C_3k\log{(h^2+k)^{-\epsilon}}}\mathbb{E}\Bigl[\textbf{1}_{\Omega_{h}^{\epsilon}}\Bigl(\sum_{n=1}^Nk\|\Pi_he_u^{n-1,h}\|^{2}_{L^2_x}\Bigr)\Bigr]\\ \nonumber&\quad
+\frac{C_4}{1-C_3k\log(h^2+k)^{-\epsilon}}\mathbb{E}\Bigl[\textbf{1}_{\Omega_{h}^{\epsilon}}\Bigl(\sum_{n=1}^Nk\|\Pi_he_u^{n-1,h}\|^2_{L^2_x}\Bigr)\Bigr]
\\ \nonumber&\leq C(h^2+k)+2\bigl({C_3\log(h^2+k)^{-\epsilon}}+{C_4}\bigr)\mathbb{E}\Bigl[\textbf{1}_{\Omega_{h}^{\epsilon}}\Bigl(\sum_{n=1}^Nk\|\Pi_he_u^{n-1,h}\|^{2}_{L^2_x}\Bigr)\Bigr].
\end{align}
Then \eqref{eq4.30} follows from an application of the discrete Gr{o}nwall inequality and the triangle inequality. The proof is complete.
\end{proof}

The second result of this section is the following error estimate for the pressure
approximation $\{r_h^n; 1\leq n\leq N\}$ and $\{p_h^n; 1\leq n\leq N\}$.

\begin{theorem}\label{thm4.5}
Let the assumptions of Theorem \ref{thm3.1} be satisfied.
Let $\{r_h^n; 1\leq n\leq {N}\}$ be the pressure approximation defined by \textbf{Algorithm 2}.
Then the following error estimate holds for $m=1,2,\cdots, N$
\begin{align}\label{eq4.50}
\mathbb{E}\Bigl[\textbf{1}_{\Omega_{h}^{\epsilon}}\Bigl(\Bigl\|k\sum^N_{n=1} \bigl(r^n- r^n_h\bigr) \Bigr\|^2_{L^2_x}\Bigr)\Bigr]\leq C(h^{2-2\epsilon}+k^{1-\epsilon}),
\end{align}
where $C$ is a positive constant independent of $h$ and $k$.
\end{theorem}

\begin{proof}
Summing \eqref{eq4.27a} over $1\leq n\leq m (\leq N)$
and subtracting the resulting equation from \eqref{eq3.14}, we have
\begin{align}\label{eq4.51}
 \bigl(e_u^{m,h},v_h\bigr) &+k\sum_{n=1}^m a\bigl(e_u^{n,h}, v_h\bigr)-k\sum_{n=1}^m d\bigl(v_h,e_r^{n,h}\bigr)\\&\nonumber
+k\sum_{n=1}^m [b\bigl(u^{n}_{h},u^{n}_{h},v_h\bigr)-b\bigl(u^{n},u^{n},v_h\bigr)]\\&\nonumber
=(e_u^0,v_h)+ \sum_{n=1}^m \bigl([{\eta}(u^{n-1})-{\eta}(u^{n-1}_h)] \Delta W_n,v_h\bigr),\quad \forall v_h\ \in \cV_h,\ a.s.
\end{align}
 Using the Poincar\'{e} inequality, the H\"{o}lder inequality and the embedding inequality,
it follows that
\begin{align}\label{eq4.52}
&d\Bigl(v_h,k\sum^m_{n=1}e_r^{n,h}\Bigr)
=\bigl(e_u^{0,h}-e_u^{n,h},v_h\bigr)-k\sum^m_{n=1}a\bigl(e_u^{n,h},v_h\bigr)\\ \nonumber
&\qquad\qquad+k\sum_{n=1}^m [b\bigl(u^{n}_{h},u^{n}_{h},v_h\bigr)-b\bigl(u^{n},u^{n},v_h\bigr)]\\ \nonumber&\qquad\qquad
+\sum^m_{n=1}\bigl([{\eta}(u^{n-1})-{\eta}(u^{n-1}_h)]\Delta W_n,v_h\bigr) \\ \nonumber&
 \leq C\Bigl(\|e_u^{0,h}\|_{L^2_x} + \|e_u^{n,h}\|_{L^2_x}
+ \sum_{n=1}^mk\|\nabla e_u^{n,h}\|_{L^2_x}\\ \nonumber&\quad
+\sum_{n=1}^mk\|\nabla e_u^{n,h}\|_{L^2_x}\|\nabla u^{n}\|_{L^2_x}+\sum_{n=1}^mk\|\nabla u^{n}_{h}\|_{L^2_x}\|\nabla e_u^{n,h}\|_{L^2_x}\\ \nonumber&\quad
+\|\sum^m_{n=1}[{\eta}(u^{n-1})-{\eta}(u^{n-1}_h)]\Delta W_n\|_{L^2_x}\Bigr)\|\nabla v\|_{L^2_x}.
\end{align}
Applying the discrete inf-sup condition \eqref{inf-sup-TH}, we obtain
\begin{align}\label{eq4.53}
\gamma\Bigl\|k\sum^m_{n=1}e_r^{n,h}\Bigr\|_{L^2_x}
& \leq C\Bigl(\|e_u^{0,h}\|_{L^2_x} + \|e_u^{n,h}\|_{L^2_x}
+ \sum_{n=1}^mk\|\nabla e_u^{n,h}\|_{L^2_x}\\ \nonumber&\quad
+\sum_{n=1}^mk\|\nabla e_u^{n,h}\|_{L^2_x}\|\nabla u^{n}\|_{L^2_x}+\sum_{n=1}^mk\|\nabla u^{n}_{h}\|_{L^2_x}\|\nabla e_u^{n,h}\|_{L^2_x}\\ \nonumber&\quad
+\Bigl\|\sum^m_{n=1}[{\eta}(u^{n-1})-{\eta}(u^{n-1}_h)]\Delta W_n\Bigr\|_{L^2_x}\Bigr).
\end{align}
With \eqref{eq4.99} and taking the expectation,  one finds that
\begin{align}\label{eq4.54}
&\quad\mathbb{E}\Bigl[\textbf{1}_{\Omega_{h}^{\epsilon}}\Bigl(\Bigl\|k\sum^m_{n=1}e_r^{n,h}\Bigr\|_{L^2_x}\Bigr)\Bigr]\\ \nonumber
& \leq C\mathbb{E}\Bigl[\textbf{1}_{\Omega_{h}^{\epsilon}}\Bigl(\|e_u^{0,h}\|_{L^2_x}\Bigr)\Bigr] + C\mathbb{E}\Bigl[\textbf{1}_{\Omega_{h}^{\epsilon}}\Bigl(\|e_u^{n,h}\|_{L^2_x}\Bigr)\Bigr]
+ C\mathbb{E}\Bigl[\textbf{1}_{\Omega_{h}^{\epsilon}}\Bigl(\sum_{n=1}^mk\|\nabla e_u^{n,h}\|_{L^2_x}\Bigr)\Bigr]\\ \nonumber&\quad
+C \mathbb{E}\Bigl[\textbf{1}_{\Omega_{h}^{\epsilon}}\Bigl(\sum_{n=1}^mk\|\nabla e_u^{n,h}\|_{L^2_x}\Bigr)\Bigl(\max_{1\leq n\leq m}\|\nabla u^{n}\|_{L^2_x}\Bigr)\Bigr] \\ \nonumber&\quad
+C \mathbb{E}\Bigl[\textbf{1}_{\Omega_{h}^{\epsilon}}\Bigl(\sum_{n=1}^mk\|\nabla e_u^{n,h}\|_{L^2_x}\Bigr)\Bigl(\max_{1\leq n\leq m}\|\nabla u^{n}_{h}\|_{L^2_x}\Bigr)\Bigr] \\ \nonumber&\quad
+C\mathbb{E}\Bigl[\textbf{1}_{\Omega_{h}^{\epsilon}}\Bigl(\Bigl\|\sum^m_{n=1}[{\eta}(u^{n-1})-{\eta}(u^{n-1}_h)]\Delta W_n\Bigr\|_{L^2_x}\Bigr)\Bigr].
\end{align}
By a standard calculation, it follows that
\begin{align}\label{eq4.55}
&\quad\mathbb{E}\Bigl[\textbf{1}_{\Omega_{h}^{\epsilon}}\Bigl(\Bigl\|k\sum^m_{n=1}e_r^{n,h}\Bigr\|_{L^2_x}\Bigr)\Bigr]\\ \nonumber
& \leq C\mathbb{E}\Bigl[\textbf{1}_{\Omega_{h}^{\epsilon}}\Bigl(\|e_u^{0,h}\|_{L^2_x}\Bigr)\Bigr] + C\mathbb{E}\Bigl[\textbf{1}_{\Omega_{h}^{\epsilon}}\Bigl(\|e_u^{n,h}\|_{L^2_x}\Bigr)\Bigr]
+ C\mathbb{E}\Bigl[\textbf{1}_{\Omega_{h}^{\epsilon}}\Bigl(\sum_{n=1}^mk\|\nabla e_u^{n,h}\|_{L^2_x}\Bigr)\Bigr]\\ \nonumber&\quad
+C\Bigl(\mathbb{E}\Bigl[\textbf{1}_{\Omega_{h}^{\epsilon}}\Bigl(\sum_{n=1}^mk\|\nabla e_u^{n,h}\|_{L^2_x}^2\Bigr)\Bigr]\Bigr)^{\frac{1}{2}}\Bigl(\mathbb{E}\Bigl[\max_{1\leq n\leq m}\|\nabla u^{n}\|_{L^2_x}^2\Bigr]\Bigr)^{\frac{1}{2}}\\ \nonumber&\quad
+C\Bigl(\mathbb{E}\Bigl[\textbf{1}_{\Omega_{h}^{\epsilon}}\Bigl(\sum_{n=1}^mk\|\nabla e_u^{n,h}\|_{L^2_x}^2\Bigr]\Bigr)^{\frac{1}{2}}\Bigl(\mathbb{E}\Bigl[\max_{1\leq n\leq m}\|\nabla u_h^{n}\|_{L^2_x}^2\Bigr)\Bigr]\Bigr)^{\frac{1}{2}}\\ \nonumber&\quad
+C\mathbb{E}\Bigl[\textbf{1}_{\Omega_{h}^{\epsilon}}\Bigl(\Bigl\|\sum^m_{n=1}[{\eta}(u^{n-1})-{\eta}(u^{n-1}_h)]\Delta W_n\Bigr\|_{L^2_x}\Bigr)\Bigr].
\end{align}
With using Lemma \ref{lem3.3} and \eqref{eq4.29a}, the last term in \eqref{eq4.55} can be bounded by \eqref{eq4.42}-\eqref{eq4.46}
which gives the desired result \eqref{eq4.50}. The proof is complete.
\end{proof}

\begin{theorem}\label{thm4.6}
Let the assumptions of Theorem \ref{thm3.1} be satisfied.
Let $\{p_h^n; 1\leq n\leq {N}\}$ be the pressure approximation defined by \textbf{Algorithm 2}.
Then the following error estimate holds for $m=1,2,\cdots, N$
\begin{align}\label{eq4.334}
\mathbb{E}\Bigl[\textbf{1}_{\Omega_{h}^{\epsilon}}\Bigl(\Bigl\|k\sum^N_{n=1} \bigl(p^n- p^n_h\bigr) \Bigr\|^2_{L^2_x}\Bigr)\Bigr]\leq C(h^{2-2\epsilon}+k^{1-\epsilon}),
\end{align}
where $C$ is a positive constant independent of $h$ and $k$.
\end{theorem}

For $\epsilon>0$, we introduce the following sample set
 \begin{align}\label{eq4.335}
 \Omega_{h,h}^{\epsilon}=\Bigl\{\omega\in \Omega\bigl|\max_{1\leq n\leq  N}\bigl(\|Au^n\|^4_{L^2_x}+\|\nabla u_h^n\|^4_{L^2_x}\bigr)\leq (h^2+k)^{-\epsilon}\Bigr\}
\end{align}
such that
 \begin{align}\label{eq4.336}
\mathbb{P}(\Omega_{h,h}^{\epsilon})\geq 1-\dfrac{\mathbb{E}\bigl[\omega\in \Omega\bigl|\max_{1\leq n\leq  N}\bigl(\|Au^n\|^4_{L^2_x}+\|\nabla u_h^n\|^4_{L^2_x}\bigr)\bigr]}{(h^2+k)^{-\epsilon}}\geq 1-\frac{C}{(h^2+k)^{-\epsilon}}.
\end{align}

The next Theorem states and proves strong optimal $H^1$ convergence for the velocity approximation.

\begin{theorem}\label{thm4.7}
Set $u^0=u_0$ and let  $\{u^n; 1\leq n\leq N\}$ and
$\{u^n_h; 1\leq n\leq N \}$ be the solutions of \textbf{ Algorithm 1} and \textbf{ Algorithm 2},
respectively.
Then, provided that $0<k<k_0$ and $0<h<h_0$ with $k_0$ and $h_0$ sufficiently small, the following error estimate holds:
\begin{align}\label{eq4.337}
& \mathbb{E}\bigl[\textbf{1}_{\Omega_{h,h}^{\epsilon}\cap\Omega_{h}^{\epsilon}}\Bigl(\max_{1\leq n\leq N}\|\nabla({u}^n-{u}^n_{h})\|^2_{L^2} 
+k\sum_{n=1}^{N}\|A({u}^n-{u}^n_{h})\|^2_{L^2}\Bigr)\Bigr] \leq C(h^{2-4\epsilon}+k^{1-2\epsilon}),
\end{align}
where $C$ is a positive constant independent of $h$ and $k$.
\end{theorem}
\begin{proof}
Taking $v_h=A_h\Pi_he_u^{n,h}\in \cV_{0h}$ and $q_h=0$ in \eqref{eq4.31}, we have
\begin{align}\label{eq4.339}
&\frac{1}{2}\bigl(\|\nabla \Pi_he_u^{n,h}\|^2_{L^2_x}-\|\nabla \Pi_he_u^{n-1,h}\|^2_{L^2_x}+\|\nabla \Pi_he_u^{n,h}-\nabla \Pi_he_u^{n-1,h}\|^2_{L^2_x}\bigr)
\\ \nonumber
& +k\,\nu\|A_h\Pi_he_u^{n,h}\|^2_{L^2_x}
=k\,\nu\bigl(A_h(u^n-u^{n}_{h}),A_h\Pi_he_u^{n,h}\bigr) +k\,b\bigl({u}_{h}^{n},{u}_{h}^{n},A_h\Pi_he_u^{n,h}\bigr)\\& \nonumber
-k\,b\bigl({u}^{n},{u}^{n},A_h\Pi_he_u^{n,h}\bigr)
+\bigl([{\eta}(u^{n-1})-{\eta}(u^{n-1}_h)]\Delta W_{n},A_h\Pi_he_u^{n,h} \bigr)\\ \nonumber&
=II_1+II_2+II_3.
\end{align}
For term $II_1$, thanks to the Young's inequality and \eqref{eq4.24}, we obtain
\begin{align}\label{eq4.340}
II_1 
\leq \frac{\nu k}{8}\|A_h\Pi_he_u^{n,h}\|^2_{L^2_x}+C{kh^2}\|\nabla^3{u}^{n}\|^2_{L^2_x}.
\end{align}
For nonlinear term $II_2$, we can decomposed as follows:
\begin{align*}
II_2&=-kb\bigl({u}^{n}-{u}_h^{n},u^{n},A_h\Pi_he_u^{n,h}\bigr)-kb\bigl(u^{n}_{h},{u}^{n}-{u}_h^{n},A_h\Pi_he_u^{n,h}\bigr)\\&
=II_{2,1}+II_{2,2}.
\end{align*}
Using the Poincar\'{e} inequality, the Young's inequality, the embedding inequality and \eqref{eq4.24}, one finds that
\begin{align}\label{eq4.341}
II_{2,1}&\leq Ck\|{u}^{n}-{u}_h^{n}\|^{\frac{1}{2}}_{L^2_x}\|\nabla({u}^{n}-{u}_h^{n})\|^{\frac{1}{2}}_{L^2_x}\|Au^{n}\|_{L^2_x}\|A_h\Pi_he_u^{n,h}\|_{L^2_x}\\&\nonumber
\leq \frac{k\nu}{8}\|A_h\Pi_he_u^{n,h}\|^2_{L^2_x}+Ck\|{u}^{n}-{u}_h^{n}\|^{2}_{L^2_x}+Ck\|\nabla({u}^{n}-{u}_h^{n})\|^2_{L^2_x}\|Au^{n}\|^4_{L^2_x},\\\label{eq4.342}
II_{2,2}&\leq k\|\nabla u^{n}_{h}\|_{L^2_x}\|\nabla ({u}^{n}-{u}_h^{n})\|^{\frac{1}{2}}_{L^2_x}\|A_h({u}^{n}-{u}_h^{n})\|^{\frac{1}{2}}_{L^2_x}\|A_h\Pi_he_u^{n,h}\|_{L^2_x}\\&\nonumber
\leq \frac{k\nu}{8}\|A_h\Pi_he_u^{n,h}\|^2_{L^2_x}+Ck\|\nabla ({u}^{n}-{u}_h^{n})\|^2_{L^2_x}\|\nabla u^{n}_{h}\|^4_{L^2_x}
+Ckh^2\|\nabla^3{u}^{n}\|^2_{L^2_x}.
\end{align}
Inserting estimates \eqref{eq4.340}--\eqref{eq4.342} into \eqref{eq4.339}, we have
\begin{align}\label{eq4.343}
&\frac{1}{2}\bigl(\|\nabla \Pi_he_u^{n,h}\|^2_{L^2_x}-\|\nabla \Pi_he_u^{n-1,h}\|^2_{L^2_x}
+\|\nabla \Pi_he_u^{n,h}-\nabla \Pi_he_u^{n-1,h}\|^2_{L^2_x}\bigr) \\\nonumber&\quad+k\,\nu\|A_h\Pi_he_u^{n,h}\|^2_{L^2_x}
\leq C{kh^2}\|\nabla^3{u}^{n}\|^2_{L^2_x}+Ck\|{u}^{n}-{u}_h^{n}\|^{2}_{L^2_x}\\ \nonumber&\qquad+Ck\|\nabla({u}^{n}-{u}_h^{n})\|^2_{L^2_x}\|Au^{n}\|^4_{L^2_x}
 +Ck\|\nabla ({u}^{n}-{u}_h^{n})\|^2_{L^2_x}\|\nabla u^{n}_{h}\|^4_{L^2_x}
\\ \nonumber&\qquad\quad+\bigl([{\eta}(u^{n-1})-{\eta}(u^{n-1}_h)]\Delta W_{n},A_h\Pi_he_u^{n,h} \bigr).
\end{align}
Taking the expectation and applying the summation operator $\sum^N_{n=1}$, one finds that
\begin{align}\label{eq4.344}
&\mathbb{E}\Big[\textbf{1}_{\Omega_{h,h}^{\epsilon}\cap\Omega_{h}^{\epsilon}}\Bigr(\frac{1}{2}\|\nabla\Pi_he_u^{n,h}\|^2_{L^2_x}+\frac{1}{2}\sum^N_{n=1}\|\nabla\Pi_he_u^{n,h}-\Pi_he_u^{n-1,h}\|^2_{L^2_x}+\sum^N_{n=1}\frac{k\nu}{2}\|A_h\Pi_he_u^{n,h}\|^2_{L^2_x}\Bigr)\Bigr]\\ \nonumber&
\leq C(h^{2-2\epsilon}+k^{1-\epsilon})(h^2+k)^{-\epsilon}+\frac{1}{2}\mathbb{E}\big[\|\nabla\Pi_he_u^{0,h}\|^2_{L^2_x}\bigr]+C{h^2}\mathbb{E}\Big[k\sum^N_{n=1}\|\nabla^3{u}^{n}\|^2_{L^2_x}\Bigr]\\ \nonumber&\qquad
+Ck\underbrace{\mathbb{E}\Big[\sum^N_{n=1}\|\nabla({u}^{n}-{u}_h^{n})\|^2_{L^2_x}\|A(u^{n}-u^{n-1})\|^4_{L^2_x}\Bigr]}_{\Lambda_5}\\ \nonumber&\qquad\qquad
+Ck\underbrace{\mathbb{E}\Big[\sum^N_{n=1}\|\nabla ({u}^{n}-{u}_h^{n})\|^2_{L^2_x}\|\nabla (u^{n}_{h}-u^{n-1}_{h})\|^4_{L^2_x}\Bigr]}_{\Lambda_6}\\ \nonumber&\qquad\qquad\qquad
+\mathbb{E}\Big[\textbf{1}_{\Omega_{h,h}^{\epsilon}\cap\Omega_{h}^{\epsilon}}\Bigr(\sum^N_{n=1}\bigl([{\eta}(u^{n-1})-{\eta}(u^{n-1}_h)]\Delta W_{n},A_h\Pi_he_u^{n,h} \bigr)\Bigr)\Bigr].
\end{align}
Now we explain how to estimate in expectation for $\Lambda_5$ and $\Lambda_6$. Using the Lemmas \ref{lem3.3}, \ref{lem3.4} and  Lemma \ref{lem4.6}, the terms $\Lambda_5$ and $\Lambda_6$ are uniformly bounded
\begin{align*}
{\Lambda_5}&
\leq  \Bigl(\mathbb{E}\max_{1\leq m\leq N}\|\nabla{u}^{m}\|^4_{L^2_x}\Bigr)^{\frac{1}{2}}\Bigl(\mathbb{E}\sum^N_{n=1}\|A(u^{n}-u^{n-1})\|^8_{L^2_x}\Bigr)^{\frac{1}{2}}\\&\nonumber\quad
 +\Bigl(\mathbb{E}\max_{1\leq m\leq N}\|\nabla{u}_h^{m}\|^4_{L^2_x}\Bigr)^{\frac{1}{2}}\Bigl(\mathbb{E}\sum^N_{n=1}\|A(u^{n}-u^{n-1})\|^8_{L^2_x}\Bigr)^{\frac{1}{2}}.
\end{align*}
and
\begin{align*}
{\Lambda_6}&
\leq \Bigl(\mathbb{E}\max_{1\leq m\leq N}\|\nabla{u}^{m}\|^4_{L^2_x}\Bigr)^{\frac{1}{2}}\Bigl(\mathbb{E}\sum^N_{n=1}\|\nabla(u^{n}_h-u^{n-1}_h)\|^8_{L^2_x}\Bigr)^{\frac{1}{2}}\\&\nonumber\quad
 +\Bigl(\mathbb{E}\max_{1\leq m\leq N}\|\nabla{u}_h^{m}\|^4_{L^2_x}\Bigr)^{\frac{1}{2}}\Bigl(\mathbb{E}\sum^N_{n=1}\|\nabla(u^{n}_h-u^{n-1}_h)\|^8_{L^2_x}\Bigr)^{\frac{1}{2}}.
\end{align*}
For term $II_3$, using the It\^{o}'s isometry and the Young's inequality, we have
\begin{align}\label{eq4.345}
&\mathbb{E}\Bigl[\Omega_{h,h}^{\epsilon}\cap\Omega_{h}^{\epsilon}\Bigr(\sum_{n=1}^{N} \bigl(\nabla[{\eta}(u^{n-1})-{\eta}(u^{n-1}_h)]\Delta W_{n},\nabla\Pi_he_u^{n,h} \bigr) \Bigr) \Bigr]\\
\nonumber
&
=\mathbb{E}\Bigl[\Omega_{h,h}^{\epsilon}\cap\Omega_{h}^{\epsilon}\Bigr(\sum_{n=1}^{N} \bigl(\nabla[{\eta}(u^{n-1})-{\eta}(u^{n-1}_h)]\Delta W_{n},\nabla\Pi_he_u^{n,h}-\nabla\Pi_he_u^{n-1,h} \bigr)\Bigr)  \Bigr]\\ \nonumber&
\leq \mathbb{E}\Bigl[\Omega_{h,h}^{\epsilon}\cap\Omega_{h}^{\epsilon}\Bigr(k\sum_{n=1}^{N} \bigl\|\nabla[{\eta}(u^{n-1})-{\eta}(u^{n-1}_h)]\Delta W_{n}\|_{L^2_x}\|\nabla(\Pi_he_u^{n,h}-\Pi_he_u^{n-1,h})\|_{L^2_x}\Bigr)  \Bigr]\\ \nonumber&
\leq \frac14 \mathbb{E}\Bigl[\Omega_{h,h}^{\epsilon}\cap\Omega_{h}^{\epsilon}\Bigr(\|\nabla(\Pi_he_u^{n,h}-\Pi_he_u^{n-1,h})\|_{L^2_x}^2\Bigr) \Bigr]+ \mathbb{E}\Bigl[\Omega_{h,h}^{\epsilon}\cap\Omega_{h}^{\epsilon}\Bigr(k\sum_{n=1}^{N} \bigl\|\nabla[{\eta}(u^{n-1})-{\eta}(u^{n-1}_h)]\|^2_{L^2_x}\Bigr) \Bigr].
\end{align}
By the definition of ${\eta}$ and using \eqref{eq2.3}, \eqref{eq2.4}, \eqref{eq2.77}, \eqref{eq3.11} and  \eqref{eq4.24}, it follows that
\begin{align}\label{eq4.346}
&\quad\bigl\|\nabla[{\eta}(u^{n-1})-{\eta}(u^{n-1}_h)]\|^2_{L^2_x}\\&\nonumber \leq C\|\nabla e_u^{n-1}\|_{L^2_x}^2+Ch^2\|{\eta}(u^{n-1})\|_{H^3_x}^2\\&\nonumber
\leq Ch^2\|\nabla^2 G({u}^{n-1})\|_{L^2}^2+C\|\nabla e_u^{n-1}\|_{L^2_x}^2\\&\nonumber
\leq Ch^2\|\nabla^2 {u}^{n-1}\|_{L^2_x}^2+Ch^4\|\nabla^3 {u}^{n-1}\|_{L^2_x}^2+C\|\nabla\Pi_he_u^{n-1,h}\|_{L^2_x}^2.
\end{align}
Combining \eqref{eq4.345}--\eqref{eq4.346} into \eqref{eq4.344}, we get
\begin{align}\label{eq4.347}
&\mathbb{E}\Big[\Omega_{h,h}^{\epsilon}\cap\Omega_{h}^{\epsilon}\Bigr(\frac{1}{2}\|\nabla\Pi_he_u^{n,h}\|^2_{L^2_x}+\frac{1}{4}\sum^N_{n=1}\|\nabla(\Pi_he_u^{n,h}-\Pi_he_u^{n-1,h})\|^2_{L^2_x}
\\ \nonumber&\quad+\sum^N_{n=1}\frac{k\nu}{2}\|A_h\Pi_he_u^{n,h}\|^2_{L^2_x}\Bigr)\Bigr]
\leq C(h^{2-4\epsilon}+k^{1-2\epsilon})+C\mathbb{E}\Bigl[\Omega_{h,h}^{\epsilon}\cap\Omega_{h}^{\epsilon}\Bigr(\sum_{n=1}^Nk\|\nabla\Pi_he_u^{n-1,h}\|_{L^2_x}^2\Bigr)\Bigr].
\end{align}
Then the \eqref{eq4.337} follows from an application of the discrete Gr{o}nwall inequality and the triangle inequality.
\end{proof}

For $\epsilon>0$, we introduce the following sample set
 \begin{align}\label{eq4.348}
 \Omega_{k,h}^{\epsilon}=\Bigl\{\omega\in \Omega\bigl|\max_{1\leq n\leq  N}\bigl(\|\nabla^2 u^n\|^4_{L^2_x}+\|\nabla u_h^n\|^4_{L^2_x}\bigr)\leq -\epsilon\log(h^4+k)\Bigr\}
\end{align}
such that
 \begin{align}\label{eq4.349}
\mathbb{P}(\Omega_{k,h}^{\epsilon})\geq 1-\dfrac{\mathbb{E}\bigl[\omega\in \Omega\bigl|\max_{1\leq n\leq  N}\bigl(\|\nabla^2 u^n\|^4_{L^2_x}+\|\nabla u_h^n\|^4_{L^2_x}\bigr)\bigr]}{-\epsilon\log(h^4+k)}\geq 1+\frac{C}{\epsilon\log(h^4+k)}.
\end{align}
For $\kappa_0>0$, the following sample set is defined as
 \begin{align}\label{eq4.3601}
 \Omega_{\kappa_0}=\Bigl\{\omega\in \Omega\bigl|\max_{1\leq n\leq  N}\|u^n-u_h^n\|^2_{L^2_x}\leq \kappa_0(h^{2-2\epsilon}+k^{1-2\epsilon})\Bigr\}.
\end{align}

The following Theorems give and derive strong optimal $L^2$ convergence of the scheme for the velocity approximation
by using the negative norm technique.

\begin{theorem}\label{thm4.8}
Set $u^0=u_0$ and let  $\{u^n; 1\leq n\leq N\}$ and
$\{u^n_h; 1\leq n\leq N \}$ be the solutions of \textbf{ Algorithm 1} and \textbf{ Algorithm 2},
respectively.
Then, provided that $0<k<k_0$ and $0<h<h_0$ with $k_0$ and $h_0$ sufficiently small, the following error estimate holds:
\begin{align}\label{eq4.350}
& \mathbb{E}\bigl[\textbf{1}_{\Omega_{k,h}^{\epsilon}\cap\Omega_{h,h}^{\epsilon}\cap\Omega_{h}^{\epsilon}\cap\Omega_{\kappa_0}}\Bigr(\max_{1\leq n\leq N}\|{u}^n-{u}^n_{h}\|^2_{-1} +k\sum_{n=1}^{N}\|{u}^n-{u}^n_{h}\|^2_{L^2_x}\Bigr)\Bigr] \leq C(\kappa_0)(h^{4-7\epsilon}+k^{1-3\epsilon}),
\end{align}
where $C(\kappa_0)$ is a positive constant independent of $h$ and $k$.
\end{theorem}
\begin{proof}
Setting $v_h=A_h^{-1}\Pi_he_u^{n,h}\in \cV_{0h}$ and $q_h=0$ in \eqref{eq4.31},
we gain
\begin{align}\label{eq4.351}
&\frac{1}{2}\bigl(\|A_h^{-\frac{1}{2}}\Pi_he_u^{n,h}\|^2_{L^2_x}-\|A_h^{-\frac{1}{2}}\Pi_he_u^{n-1,h}\|^2_{L^2_x}+\|A_h^{-\frac{1}{2}}\Pi_he_u^{n,h}-A_h^{-\frac{1}{2}}\Pi_he_u^{n-1,h}\|^2_{L^2_x}\bigr) \\ \nonumber
&\quad+k\,\nu\| \Pi_he_u^{n,h}\|^2_{L^2_x}
=k\,a\bigl(u^n-u^{n}_{h},A_h^{-1}\Pi_he_u^{n,h}\bigr)+k\,b\bigl({u}_{h}^{n-1},{u}_{h}^{n},A_h^{-1}\Pi_he_u^{n,h}\bigr)\\ \nonumber
&\quad
-k\,b\bigl({u}^{n-1},{u}^{n},A_h^{-1}\Pi_he_u^{n,h}\bigr)+\bigl([{\eta}(u^{n-1})-{\eta}(u^{n-1}_h)]\Delta W_{n},A_h^{-1}\Pi_he_u^{n,h} \bigr)\\ \nonumber&
=III_1+III_2+III_3.
\end{align}
For term $III_1$, thanks to the Young's inequality, \eqref{eq4.23} and \eqref{eq4.24}, we obtain
\begin{align}\label{eq4.352}
III_1 
\leq \frac{\nu k}{8}\|\Pi_he_u^{n,h}\|^2_{L^2_x}+C{kh^4}\|\nabla^2{u}^{n}\|^2_{L^2_x}.
\end{align}
For nonlinear term $III_2$, we can decomposed as follows:
\begin{align*}
III_2&=-kb\bigl({u}^{n}-{u}_h^{n},u^{n},A_h^{-1}\Pi_he_u^{n,h}\bigr)-kb\bigl(u^{n}_{h},{u}^{n}-{u}_h^{n}, A_h^{-1}\Pi_he_u^{n,h}\bigr)\\&
=III_{2,1}+III_{2,2}.
\end{align*}
Using the Poincar\'{e} inequality, the Young's inequality and the embedding inequality, one finds that
\begin{align}\label{eq4.353}
III_{2,1}&\leq \frac{k\nu}{8}\|\Pi_he_u^{n,h}\|^2_{L^2_x}+Ck\|\nabla u^n\|^{2}_{L^2_x}\|{u}^{n}-\Pi_h{u}^{n}\|^{2}_{L^2_x}\\\nonumber
&\quad +Ck\|\nabla^2u^n_h\|^2_{L^2_x}\|\Pi_he_u^{n,h}\|^2_{-1}+Ck\|\Pi_he_u^{n,h}\|^2_{L^2_x}\|\nabla(u^n-u^n_h)\|^2_{L^2_x},\\\label{eq4.354}
III_{2,2}&\leq \frac{1}{4}\|\Pi_he_u^{n,h}-\Pi_he_u^{n-1,h}\|^2_{-1}+\frac{k\nu}{8}\|\Pi_he_u^{n,h}\|^2_{L^2_x}\\&\quad\nonumber
+Ck\|u^n-\Pi_hu^n\|^{2}_{L^2_x}+Ck\|\nabla^2u^n_h\|^2_{L^2_x}\|\Pi_he_u^{n-1,h}\|^2_{-1}\\\nonumber
&\quad +Ck^2\|\nabla u^{n}_h\|^2_{L^2_x}\|\nabla(u^n-u^n_h)\|^2_{L^2_x}.
\end{align}
Inserting estimates \eqref{eq4.352}--\eqref{eq4.354} into \eqref{eq4.351}, we have
\begin{align}\label{eq4.355}
&\frac{1}{2}\bigl(\|\Pi_he_u^{n,h}\|^2_{-1}-\|\Pi_he_u^{n-1,h}\|^2_{-1}+\|\Pi_he_u^{n,h}-\Pi_he_u^{n-1,h}\|^2_{-1}\bigr)\\ \nonumber&\quad
+\frac{k\nu}{2}\|\Pi_he_u^{n,h}\|^2_{L^2_x}
\leq C{kh^4}\|\nabla^2{u}^{n}\|^2_{L^2_x}+Ck\|\nabla^2u^n_h\|^2_{L^2_x}\|\Pi_he_u^{n,h}\|^2_{-1}\\ \nonumber&\qquad
+Ck\|\Pi_he_u^{n,h}\|^2_{L^2_x}\|\nabla(u^n-u^n_h)\|^2_{L^2_x}
+Ck^2\|\nabla u^{n}_h\|^2_{L^2_x}\|\nabla(u^n-u^n_h)\|^2_{L^2_x}\\ \nonumber &\qquad\quad
+\bigl([{\eta}(u^{n-1})-{\eta}(u^{n-1}_h)]\Delta W_{n},A^{-1}_h\Pi_he_u^{n,h} \bigr)
\end{align}
Taking the expectation and applying the summation operator $\sum^N_{n=1}$, one finds that
\begin{align}\label{eq4.356}
&\mathbb{E}\Big[\textbf{1}_{\Omega_{k,h}^{\epsilon}\cap\Omega_{h,h}^{\epsilon}\cap\Omega_{h}^{\epsilon}\cap\Omega_{\kappa_0}}\Bigr(\frac{1}{2}\|\Pi_he_u^{n,h}\|^2_{-1}+\frac{1}{2}\sum^N_{n=1}\|\Pi_he_u^{n,h}-\Pi_he_u^{n-1,h}\|^2_{-1}
+\sum^N_{n=1}\frac{k\nu}{4}\|\Pi_he_u^{n,h}\|^2_{L^2_x}\Bigr)\Bigr]\\ \nonumber&
\leq \frac{1}{2}\mathbb{E}\big[\|\Pi_he_u^{0,h}\|^2_{-1}\bigr]+C{h^4}\mathbb{E}\Big[k\sum^N_{n=1}\|\nabla{u}^{n}\|^2_{L^2_x}\Bigr]\\ \nonumber&\quad
+C\log(h^4+k)^{-\epsilon}\mathbb{E}\Big[\textbf{1}_{\Omega_{k,h}^{\epsilon}\cap\Omega_{h,h}^{\epsilon}\cap\Omega_{h}^{\epsilon}\cap\Omega_{\kappa_0}}\Bigr(
k\sum^N_{n=1}\|\Pi_he_u^{n,h}\|^2_{-1}\Bigr)\Bigr]\\ \nonumber&\quad
+Ck(h^{2-4\epsilon}+k^{1-2\epsilon})+Ck\underbrace{\mathbb{E}\Big[\sum^N_{n=1}\|\nabla^2(u^n_h-u^{n-1}_h)\|^2_{L^2_x}\|u^n-u^n_h\|^2_{-1}\Bigr]}_{\Lambda_7}\\ \nonumber&\quad
+C(\kappa_0)(h^{2-4\epsilon}+k^{1-2\epsilon})(h^{2-2\epsilon}+k^{1-\epsilon})+Ck^2\underbrace{\mathbb{E}\Big[\sum^N_{n=1}\|\nabla (u^n_h-u_h^{n-1})\|^{2}_{L^2_x}\|\nabla(u^n-u^n_h)\|^{2}_{L^2_x}\Bigr]}_{\Lambda_8}\\ \nonumber&\quad
+\mathbb{E}\Big[\textbf{1}_{\Omega_{k,h}^{\epsilon}\cap\Omega_{h,h}^{\epsilon}\cap\Omega_{h}^{\epsilon}\cap\Omega_{\kappa_0}}\Bigr(\sum^N_{n=1}\bigl(A^{-\frac12}_h[{\eta}(u^{n-1})-{\eta}(u^{n-1}_h)]\Delta W_{n},A^{-\frac12}_h\Pi_he_u^{n,h} \bigr)\Bigr)\Bigr].
\end{align}
Now we explain how to estimate in expectation for $\Lambda_7$ and $\Lambda_8$. Making use of the Lemma \ref{lem3.4} and Lemma \ref{lem4.5}, the terms $\Lambda_7$ and $\Lambda_8$ are uniformly bounded
\begin{align*}
{\Lambda_7}&
\leq \Bigl(\mathbb{E}\max_{1\leq m\leq N}\|{u}^{m}\|^4_{-1}\Bigr)^{\frac{1}{2}}\Bigl(\mathbb{E}\sum^N_{n=1}\|\nabla^2(u^{n}_h-u_h^{n-1})\|^4_{L^2_x}\Bigr)^{\frac{1}{2}}\\&\nonumber\quad
+\Bigl(\mathbb{E}\max_{1\leq m\leq N}\|{u}_h^{m}\|^4_{-1}\Bigr)^{\frac{1}{2}}\Bigl(\mathbb{E}\sum^N_{n=1}\|\nabla^2(u_h^{n}-u_h^{n-1})\|^4_{L^2_x}\Bigr)^{\frac{1}{2}},\\\nonumber
{\Lambda_8}&
\leq \Bigl(\mathbb{E}\max_{1\leq m\leq N}\|\nabla{u}^{m}\|^4_{L^2_x}\Bigr)^{\frac{1}{2}}\Bigl(\mathbb{E}\sum^N_{n=1}\|\nabla(u^{n}-u^{n-1})\|^4_{L^2_x}\Bigr)^{\frac{1}{2}}\\&\nonumber\quad
+\Bigl(\mathbb{E}\max_{1\leq m\leq N}\|\nabla{u}_h^{m}\|^4_{L^2_x}\Bigr)^{\frac{1}{2}}\Bigl(\mathbb{E}\sum^N_{n=1}\|\nabla(u^{n}-u^{n-1})\|^4_{L^2_x}\Bigr)^{\frac{1}{2}}.
\end{align*}
For term $III_3$, using the It\^{o}'s isometry and the Young's inequality, we have
\begin{align}\label{eq4.357}
&\quad\mathbb{E}\Bigl[\textbf{1}_{\Omega_{k,h}^{\epsilon}\cap\Omega_{h,h}^{\epsilon}\cap\Omega_{h}^{\epsilon}\cap\Omega_{\kappa_0}}\Bigr(\sum_{n=1}^{N} \bigl(A^{-\frac12}_h[{\eta}(u^{n-1})-{\eta}(u^{n-1}_h)]\Delta W_{n},\Pi_hA^{-\frac12}_he_u^{n,h} \bigr) \Bigr)\Bigr]\\
\nonumber
&
=\mathbb{E}\Bigl[\textbf{1}_{\Omega_{k,h}^{\epsilon}\cap\Omega_{h,h}^{\epsilon}\cap\Omega_{h}^{\epsilon}\cap\Omega_{\kappa_0}}\Bigr(\sum_{n=1}^{N} \bigl(A^{-\frac12}_h[{\eta}(u^{n-1})-{\eta}(u^{n-1}_h)]\Delta W_{n},\Pi_hA^{-\frac12}_he_u^{n,h}-\Pi_hA^{-\frac12}_he_u^{n-1,h} \bigr) \Bigr)\Bigr]\\ \nonumber&
\leq \mathbb{E}\Bigl[\textbf{1}_{\Omega_{k,h}^{\epsilon}\cap\Omega_{h,h}^{\epsilon}\cap\Omega_{h}^{\epsilon}\cap\Omega_{\kappa_0}}\Bigr(k\sum_{n=1}^{N} \bigl\|A^{-\frac12}_h[{\eta}(u^{n-1})-{\eta}(u^{n-1}_h)]\Delta W_{n}\|_{L^2_x}\|\Pi_hA^{-\frac12}_he_u^{n,h}-\Pi_hA^{-\frac12}_he_u^{n-1,h}\|_{L^2_x}\Bigr) \Bigr]\\ \nonumber&
\leq \frac14 \mathbb{E}\Bigl[\textbf{1}_{\Omega_{k,h}^{\epsilon}\cap\Omega_{h,h}^{\epsilon}\cap\Omega_{h}^{\epsilon}\cap\Omega_{\kappa_0}}\Bigr(\|\Pi_he_u^{n,h}-\Pi_he_u^{n-1,h}\|^2_{-1}\Bigr)\Bigr]
\\ \nonumber&\quad+ \mathbb{E}\Bigl[\textbf{1}_{\Omega_{k,h}^{\epsilon}\cap\Omega_{h,h}^{\epsilon}\cap\Omega_{h}^{\epsilon}\cap\Omega_{\kappa_0}}\Bigr(k\sum_{n=1}^{N} \bigl\|{\eta}(u^{n-1})-{\eta}(u^{n-1}_h)\|^2_{-1}\Bigr) \Bigr].
\end{align}
By the definition of ${\eta}$ and using \eqref{eq2.3}, \eqref{eq2.4}, \eqref{eq2.77}, \eqref{eq3.11} and  \eqref{eq4.24}, one finds that
\begin{align}\label{eq4.358}
&\quad\|{\eta}(u^{n-1})-{\eta}(u_h^{n-1})\|^2_{-1}\\&\nonumber \leq C\|e_u^{n-1}\|_{-1}^2+Ch^4\|{\eta}(u^{n-1})\|_{H^1_x}^2\\&\nonumber
\leq Ch^4\|\nabla\cdot G({u}^{n-1})\|_{L^2}^2+C\|e_u^{n-1}\|_{-1}^2\\&\nonumber
\leq Ch^4\|\nabla {u}^{n-1}\|_{L^2_x}^2+Ch^4\|\nabla^2 {u}^{n-1}\|_{L^2_x}^2+C\|\Pi_he_u^{n-1,h}\|_{-1}^2.
\end{align}
Combining \eqref{eq4.357}--\eqref{eq4.358} into \eqref{eq4.356}, we get
\begin{align}\label{eq4.359}
&\mathbb{E}\Big[\textbf{1}_{\Omega_{k,h}^{\epsilon}\cap\Omega_{h,h}^{\epsilon}\cap\Omega_{h}^{\epsilon}\cap\Omega_{\kappa_0}}\Bigr(\frac{1}{2}\|\Pi_he_u^{n,h}\|^2_{-1}+\frac{1}{2}\sum^N_{n=1}\|\Pi_he_u^{n,h}-\Pi_he_u^{n-1,h}\|^2_{-1}
+\sum^N_{n=1}\frac{k\nu}{4}\|\Pi_he_u^{n,h}\|^2_{L^2_x}\Bigr)\Bigr]\\ \nonumber&
\leq C(\kappa_0)(h^{4-6\epsilon}+k^{1-2\epsilon})
+C\log(h^4+k)^{-\epsilon}\mathbb{E}\Big[\textbf{1}_{\Omega_{k,h}^{\epsilon}\cap\Omega_{h,h}^{\epsilon}\cap\Omega_{h}^{\epsilon}\cap\Omega_{\kappa_0}}\Bigr(k\sum^N_{n=1}\|\Pi_he_u^{n,h}\|^2_{-1}\Bigr)\Bigr]\\ \nonumber&\quad
+\mathbb{E}\Big[\textbf{1}_{\Omega_{k,h}^{\epsilon}\cap\Omega_{h,h}^{\epsilon}\cap\Omega_{h}^{\epsilon}\cap\Omega_{\kappa_0}}\Bigr(k\sum^N_{n=1}\|\Pi_he_u^{n-1,h}\|^2_{-1}\Bigr)\Bigr].
\end{align}
Using the similar line in the proof of Theorem \ref{thm4.4},
with applying the discrete Gr{o}nwall inequality and the triangle inequality, the result \eqref{eq4.350} holds. The proof is complete.
\end{proof}

For $\kappa>0$, we introduce  the following sample set
 \begin{align}\label{eq4.360}
 \Omega_{\kappa}=\Bigl\{\omega\in \Omega\bigl|\max_{1\leq n\leq  N}\|\nabla (u^n-u_h^n)\|^2_{L^2_x}\leq \kappa(h^{2-4\epsilon}+k^{1-2\epsilon})\Bigr\}.
\end{align}

\begin{theorem}\label{thm4.9}
Set $u^0=u_0$ and let  $\{u^n; 1\leq n\leq N\}$ and
$\{u^n_h; 1\leq n\leq N \}$ be the solutions of \textbf{ Algorithm 1} and \textbf{ Algorithm 2},
respectively.
Then, provided that $0<k<k_0$ and $0<h<h_0$ with $k_0$ and $h_0$ sufficiently small, the following error estimate holds:
\begin{align}\label{eq4.461}
& \mathbb{E}\bigl[\textbf{1}_{\Omega_{k,h}^{\epsilon}\cap\Omega_{h,h}^{\epsilon}\cap\Omega_{h}^{\epsilon}\cap\Omega_{\kappa}\cap\Omega_{\kappa_0}}\bigr(\|{u}^n-{u}^n_{h}\|^2_{L^2}\bigr)\bigr] \leq C(\kappa_0,\kappa)(h^{4-8\epsilon}+k^{1-4\epsilon}),
\end{align}
where $C(\kappa_0,\kappa)$ is a positive constant independent of $h$ and $k$.
\end{theorem}
\begin{proof}
Taking $v_h=\Pi_he_u^{n,h}\in \cV_{0h}$ and $q_h=0$, we have
\begin{align}\label{eq4.361}
\frac{1}{2}\bigl(\|\Pi_he_u^{n,h}\|^2_{L^2_x}&-\|\Pi_he_u^{n-1,h}\|^2_{L^2_x}+\|\Pi_he_u^{n,h}-\Pi_he_u^{n-1,h}\|^2_{L^2_x}\bigr) \\ \nonumber
&\quad+k\,\nu\|\nabla \Pi_he_u^{n,h}\|^2_{L^2_x}
=k\,a\bigl(u^n-u^{n}_{h},\Pi_he_u^{n,h}\bigr) \\ \nonumber
&\quad+k\,b\bigl({u}_{h}^{n-1},{u}_{h}^{n},\Pi_he_u^{n,h}\bigr)-k\,b\bigl({u}^{n-1},{u}^{n},\Pi_he_u^{n,h}\bigr)\\ \nonumber
&\quad+\bigl([{\eta}(u^{n-1})-{\eta}(u^{n-1}_h)]\Delta W_{n},\Pi_he_u^{n,h} \bigr)\\ \nonumber&
=IV_1+IV_2+IV_3.
\end{align}
For term $IV_1$, thanks to the Young's inequality, \eqref{eq4.23} and \eqref{eq4.24}, we obtain
\begin{align}\label{eq4.362}
IV_1 &
\leq \frac{\nu k}{4}\|\nabla \Pi_he_u^{n,h}\|^2_{L^2_x}+C{kh^4}\|\nabla^3{u}^{n}\|^2_{L^2_x}.
\end{align}
For nonlinear term $IV_2$, using the Poincar\'{e} inequality, the Young's inequality and the embedding inequality, one finds that
\begin{align}\label{eq4.3622}
IV_2&=-kb\bigl({u}^{n}-{u}^{n}_h,u^{n},\Pi_he_u^{n,h}\bigr)-kb\bigl(u^{n}_h,{u}^{n}-{u}^{n}_h,\Pi_he_u^{n,h}\bigr)\\&\nonumber
\leq \frac{k\nu}{4}\|\nabla\Pi_he_u^{n,h}\|^2_{L^2_x}+Ck\|\nabla^2 u^{n}\|^{2}_{L^2_x}\|{u}^{n}-{u}^{n}_h\|^2_{L^2_x}\\ \nonumber&\quad
+Ck\|\nabla({u}^{n}-{u}^{n}_h)\|^4_{L^2_x}.
\end{align}
Inserting estimates \eqref{eq4.362}--\eqref{eq4.3622} into \eqref{eq4.361}, we have
\begin{align}\label{eq4.363}
&\frac{1}{2}\bigl(\|\Pi_he_u^{n,h}\|^2_{L^2_x}-\|\Pi_he_u^{n-1,h}\|^2_{L^2_x}+\|\Pi_he_u^{n,h}-\Pi_he_u^{n-1,h}\|^2_{L^2_x}\bigr)\\ \nonumber&\quad
+\frac{k\nu}{2}\|\nabla\Pi_he_u^{n,h}\|^2_{L^2_x}
\leq C{kh^4}\|\nabla^3{u}^{n}\|^2_{L^2_x}+Ck\|\nabla^2 u^{n}\|^{2}_{L^2_x}\|{u}^{n}-{u}^{n}_h\|^2_{L^2_x}\\ \nonumber&\qquad
+Ck\|\nabla({u}^{n}-{u}^{n}_h)\|^4_{L^2_x}
+\bigl([{\eta}(u^{n-1})-{\eta}(u^{n-1}_h)]\Delta W_{n},\Pi_he_u^{n,h} \bigr)
\end{align}
Taking the expectation and applying the summation operator $\sum^N_{n=1}$, it follows that
\begin{align}\label{eq4.364}
&\mathbb{E}\Big[\textbf{1}_{\Omega_{k,h}^{\epsilon}\cap\Omega_{h,h}^{\epsilon}\cap\Omega_{h}^{\epsilon}\cap\Omega_{\kappa}\cap\Omega_{\kappa_0}}\Bigr(\frac{1}{2}\|\Pi_he_u^{n,h}\|^2_{L^2_x}+\frac{1}{2}\sum^N_{n=1}\|\Pi_he_u^{n,h}
-\Pi_he_u^{n-1,h}\|^2_{L^2_x}+\sum^N_{n=1}\frac{k\nu}{2}\|\nabla \Pi_he_u^{n,h}\|^2_{L^2_x}\Bigr)\Bigr]\\ \nonumber&
\leq \frac{1}{2}\mathbb{E}\big[\|\Pi_he_u^{0,h}\|^2_{L^2_x}\bigr]+C{h^4}\mathbb{E}\Big[k\sum^N_{n=1}\|\nabla^3{u}^{n}\|^2_{L^2_x}\Bigr]\\ \nonumber&\quad
+C\mathbb{E}\Big[\textbf{1}_{\Omega_{k,h}^{\epsilon}\cap\Omega_{h,h}^{\epsilon}\cap\Omega_{h}^{\epsilon}\cap\Omega_{\kappa}\cap\Omega_{\kappa_0}}\Bigr(k\sum^N_{n=1}\|\nabla^2 u^{n}\|^{2}_{L^2_x}\|{u}^{n}-{u}^{n}_h\|^2_{L^2_x}\Bigr)\Bigr] \\ \nonumber&\quad
+C(\kappa_0,\kappa)(h^{2-4\epsilon}
+k^{1-2\epsilon})\mathbb{E}\Big[\textbf{1}_{\Omega_{k,h}^{\epsilon}\cap\Omega_{h,h}^{\epsilon}\cap\Omega_{h}^{\epsilon}\cap\Omega_{\kappa}\cap\Omega_{\kappa_0}}\Bigr(k\sum^N_{n=1}\|\nabla({u}^{n}-{u}^{n}_h)\|^2_{L^2_x}\Bigr)\Bigr]\\ \nonumber&\quad
+\mathbb{E}\Big[\textbf{1}_{\Omega_{k,h}^{\epsilon}\cap\Omega_{h,h}^{\epsilon}\cap\Omega_{h}^{\epsilon}\cap\Omega_{\kappa}\cap\Omega_{\kappa_0}}\Bigr(\sum^N_{n=1}\bigl([{\eta}(u^{n-1})-{\eta}(u^{n-1}_h)]\Delta W_{n},\Pi_he_u^{n,h} \bigr)\Bigr)\Bigr]\\ \nonumber&
\leq \frac{1}{2}\mathbb{E}\big[\|\Pi_he_u^{0,h}\|^2_{L^2_x}\bigr]+C(\kappa_0,\kappa)(h^{4-8\epsilon}+k^{1-4\epsilon})
+Ck\underbrace{\mathbb{E}\Big[\sum^N_{n=1}\|\nabla^2(u^{n}-u^{n-1})\|^{2}_{L^2_x}\|{u}^{n}-{u}^{n}_h\|^2_{L^2_x}\Bigr]}_{\Lambda_9} \\ \nonumber&\quad
+\mathbb{E}\Big[\textbf{1}_{\Omega_{k,h}^{\epsilon}\cap\Omega_{h,h}^{\epsilon}\cap\Omega_{h}^{\epsilon}\cap\Omega_{\kappa}\cap\Omega_{\kappa_0}}\Bigr(\sum^N_{n=1}\bigl([{\eta}(u^{n-1})-{\eta}(u^{n-1}_h)]\Delta W_{n},\Pi_he_u^{n,h} \bigr)\Bigr)\Bigr].
\end{align}
Using the Lemma \ref{lem3.4} and Lemma \ref{lem4.5}, the term $\Lambda_9$ is uniformly bounded. Then we obtain
\begin{align}\label{eq4.365}
&\mathbb{E}\Big[\textbf{1}_{\Omega_{k,h}^{\epsilon}\cap\Omega_{h,h}^{\epsilon}\cap\Omega_{h}^{\epsilon}\cap\Omega_{\kappa}\cap\Omega_{\kappa_0}}\Bigr(\frac{1}{2}\|\Pi_he_u^{n,h}\|^2_{L^2_x}+\frac{1}{4}\sum^N_{n=1}\|\Pi_he_u^{n,h}-\Pi_he_u^{n-1,h}\|^2_{L^2_x}
+\sum^N_{n=1}\frac{k\nu}{2}\|\nabla e_u^{n,h}\|^2_{L^2_x}\Bigr)\Bigr]\\ \nonumber&
\leq C(\kappa_0,\kappa)(h^{4-8\epsilon}+k^{1-4\epsilon})
+C\mathbb{E}\Bigl[\textbf{1}_{\Omega_{k,h}^{\epsilon}\cap\Omega_{h,h}^{\epsilon}\cap\Omega_{h}^{\epsilon}\cap\Omega_{\kappa}\cap\Omega_{\kappa_0}}\Bigr(\sum_{n=1}^Nk\|\Pi_he_u^{n-1,h}\|_{L^2_x}^2\Bigr)\Bigr].
\end{align}
By applying the discrete Gr{o}nwall inequality and the triangle inequality, the result \eqref{eq4.461} holds. The proof is complete.
\end{proof}

Theorems \ref{thm3.1}, \ref{thm3.2}, \ref{thm3.3}, \ref{thm4.5}, \ref{thm4.6}, \ref{thm4.7} and Theorem \ref{thm4.9} and the triangle inequality
infer the global error estimates, which are the main results of this paper.

\begin{theorem}\label{thm4.10}
Under the assumptions of Theorems \ref{thm3.1}, \ref{thm3.2}, \ref{thm3.3}, \ref{thm4.5}, \ref{thm4.6}, \ref{thm4.7} and Theorem \ref{thm4.9},
there hold the following error estimates:
\begin{align}\label{eq4.366}
\qquad\mathbb{E}\Bigl[\textbf{1}_{\Omega_{k,h}^{\epsilon}\cap\Omega_{h,h}^{\epsilon}\cap\Omega_{h}^{\epsilon}\cap\Omega_{\tau}^{\epsilon}}\Bigl( \|\nabla (u(t_n)-u^n_h)\|^2_{L^2_x}
\Bigr)\Bigr]&\leq C\bigl(k^{2\alpha-2\epsilon}+h^{2-4\epsilon}\bigr), \\\label{eq4.367}
\mathbb{E}\Bigl[\textbf{1}_{\Omega_{k}^{\epsilon}\cap\Omega_{h}^{\epsilon}}\Bigl(\Bigl\|\int_0^{t_m}r(s)\,ds -k\sum^m_{n=1}r^n_h \Bigr\|^2_{L^2_x}&+\Bigl\|\int_0^{t_m} p(s)\,ds -k\sum^m_{n=1}p^n_h \Bigr\|^2_{L^2_x}\Bigr)\Bigr]
\\\nonumber &\leq C\bigl(k^{2\alpha-\epsilon}+h^{2-2\epsilon}\bigr),\\\label{eq4.368}
\mathbb{E}\Bigl[\textbf{1}_{\Omega_{k,h}^{\epsilon}\cap\Omega_{h,h}^{\epsilon}\cap\Omega_{h}^{\epsilon}\cap\Omega_{k}^{\epsilon}\cap\Omega_{\kappa_0}\cap\Omega_{\kappa}}\Bigl( \|u(t_n)-u^n_h\|^2_{L^2_x}
\Bigr)\Bigr]&\leq C(\kappa_0,\kappa)\bigl(k^{2\alpha-4\epsilon}+h^{4-8\epsilon}\bigr),
\end{align}
where $C,C(\kappa_0,\kappa)$ are two positive constants independent of $h$ and $k$.
\end{theorem}

\begin{remark}
The crucial point which makes the error analysis interesting and distinct from the deterministic
case is the low regularity in time. As far as the spatial
regularity  is concerned, we can obtain similar optimal error estimates to the deterministic case.
From the numerical results of Section 5,  the estimates \eqref{eq4.366}--\eqref{eq4.368} are optimal order.
 \end{remark}

\section{Numerical results}\label{sec-5}
In this section, we give some $2$D numerical results to confirm the theoretical
error estimates of our \textbf{Algorithm 2}.
We set $D=(0,L)^2$ with $L=1$, a deterministic constant force term, the initial condition $u_0=(0,0)$
and $G(u(t))= \bigl((u_1(t)^2+1)^{\frac{1}{2}},(u_2(t)^2+1)^{\frac{1}{2}}\bigr)$. 
The $W$ in \eqref{eq1.1} is chosen as a finite-dimensional ${\mathbb{R}}^M$-Wiener
process such that
\[
W(t_{n})-W(t_{n-1})=\sum_{j=1}^M \sum_{k=1}^M \lambda^{\frac{1}{2}}_{j,k}g_{j,k}\xi^n_{j,k},
\]
where $\lambda_{j,k}=\frac{1}{j^2+k^2}$, $\xi^n_{j,k}\sim N(0,1)$ and
$g_{j,k}(x,y)=5\sin(j\pi x)\sin(k\pi y)$.

We take the following parameters: $M=10$, $\nu=1$ and $T=1$.
 The Monte Carlo method
with $N_p = 1200$ realizations is utilized to compute the expectation.
Since the exact solution of the problem \eqref{eq1.1} is unknown,
 we denote the time/spatial errors of the numerical
solutions by
\begin{align*}
EAu^N_k &:=\Bigl(\mathbb{E}\Bigl[\|u_{h,k}^N-u_{h,k/2}^N\|^2_{L^2}\Bigr]\Bigr)^{\frac12},   &&EAu^N_h:=\Bigl(\mathbb{E}\Bigl[\|u_{h,k}^N-u_{h/2,k}^N\|^2_{L^2}\Bigr]\Bigr)^{\frac12},\\
EBu^N_k &:=\Bigl(\mathbb{E}\Bigl[\|\nabla (u_{h,k}^N-u_{h,k/2}^N)\|^2_{L^2})\Bigr]\Bigr)^{\frac12},    && EBu^N_h:=\Bigl(\mathbb{E}\Bigl[\|\nabla (u_{h,k}^N-u_{h/2,k}^N)\|^2_{L^2})\Bigr]\Bigr)^{\frac12},\\
Ep^N_k:&=\Bigl(\mathbb{E}\Bigl[\Bigl\|{k}\sum_{n=1}^{\frac{T}{k}}p_{h}^n-\frac{k}{2}\sum_{n=1}^{\frac{2T}{k}}p_{h}^n\Bigr\|^2_{L^2}\Bigr]\Bigr)^{\frac12},  &&Ep^N_h:=\Bigl(\mathbb{E}\Bigl[\Bigl\|{k}\sum_{n=1}^{\frac{T}{k}}p_{h}^n-{k}\sum_{n=1}^{\frac{T}{k}}p_{h/2}^n\Bigr\|^2_{L^2}\Bigr]\Bigr)^{\frac12}.
\end{align*}
where $u_{h,k}^N$ is the one path simulation at $t^N = T$ computed by using space mesh size $h$ and time mesh size $k$.

Figure 1 shows the errors of the time discretizations 
of the velocity and the pressure
using different time mesh sizes. It is evident that the numerical results validate
the half order for the time discretization as theoretical error estimates.
Figure 2 displays the errors of the spatial discretizations 
of the velocity and the pressure
using different space mesh sizes. It is easy to see that the numerical
results check the first/second order for the spatial discretization as proved in Theorem \ref{thm4.10}.
\begin{figure}[th]
\centerline{
\includegraphics[scale=0.28]{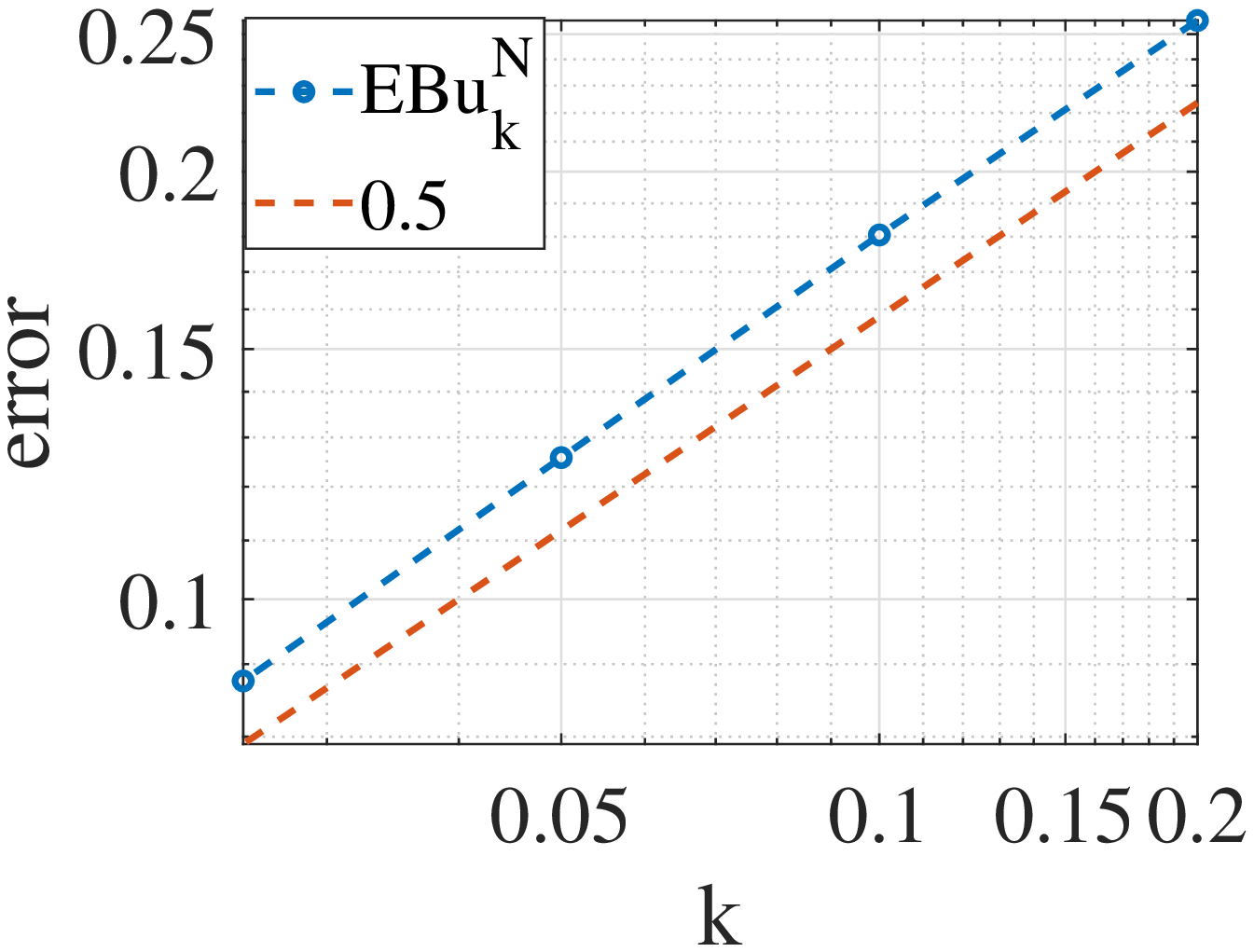}
\includegraphics[scale=0.28]{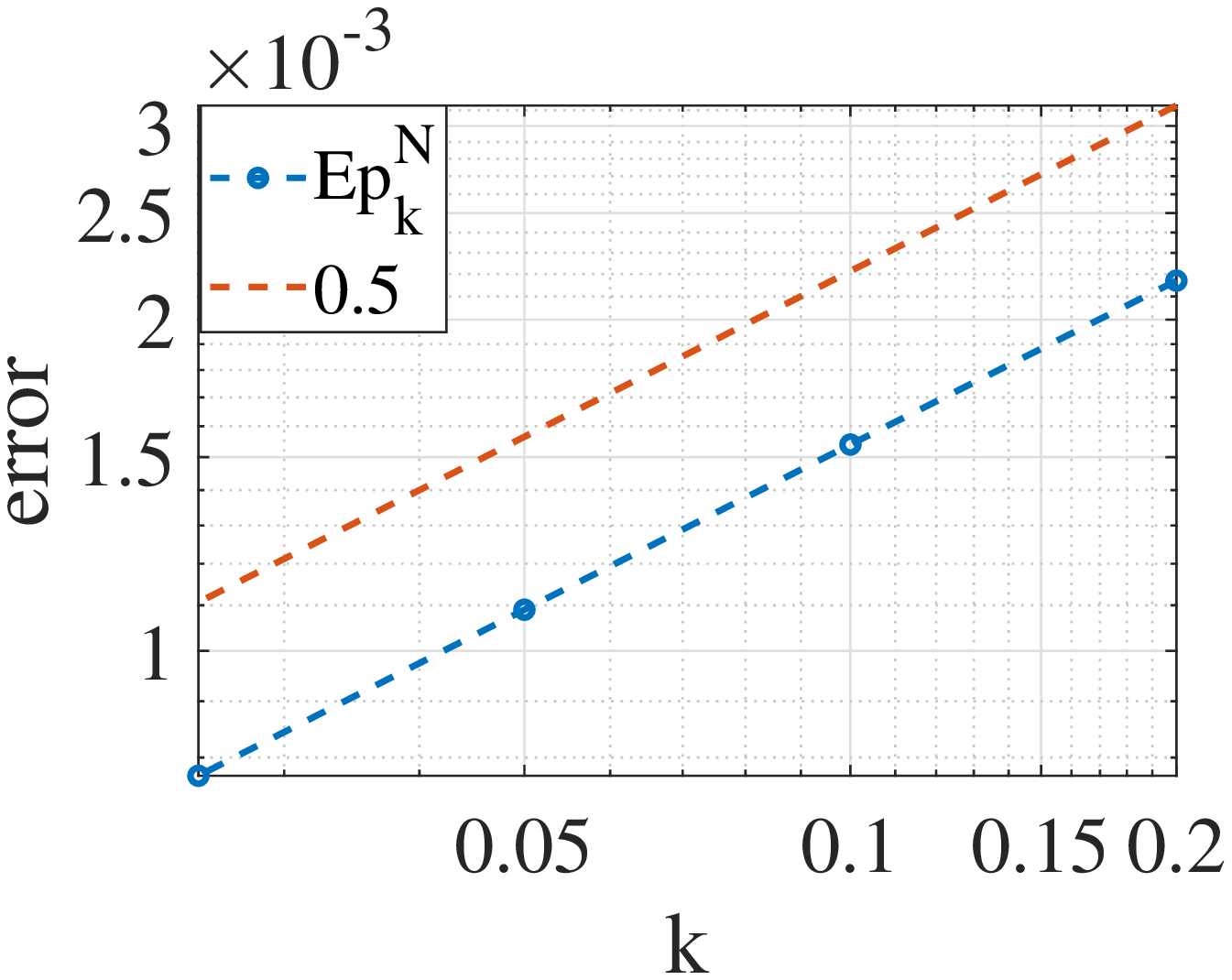}
\includegraphics[scale=0.28]{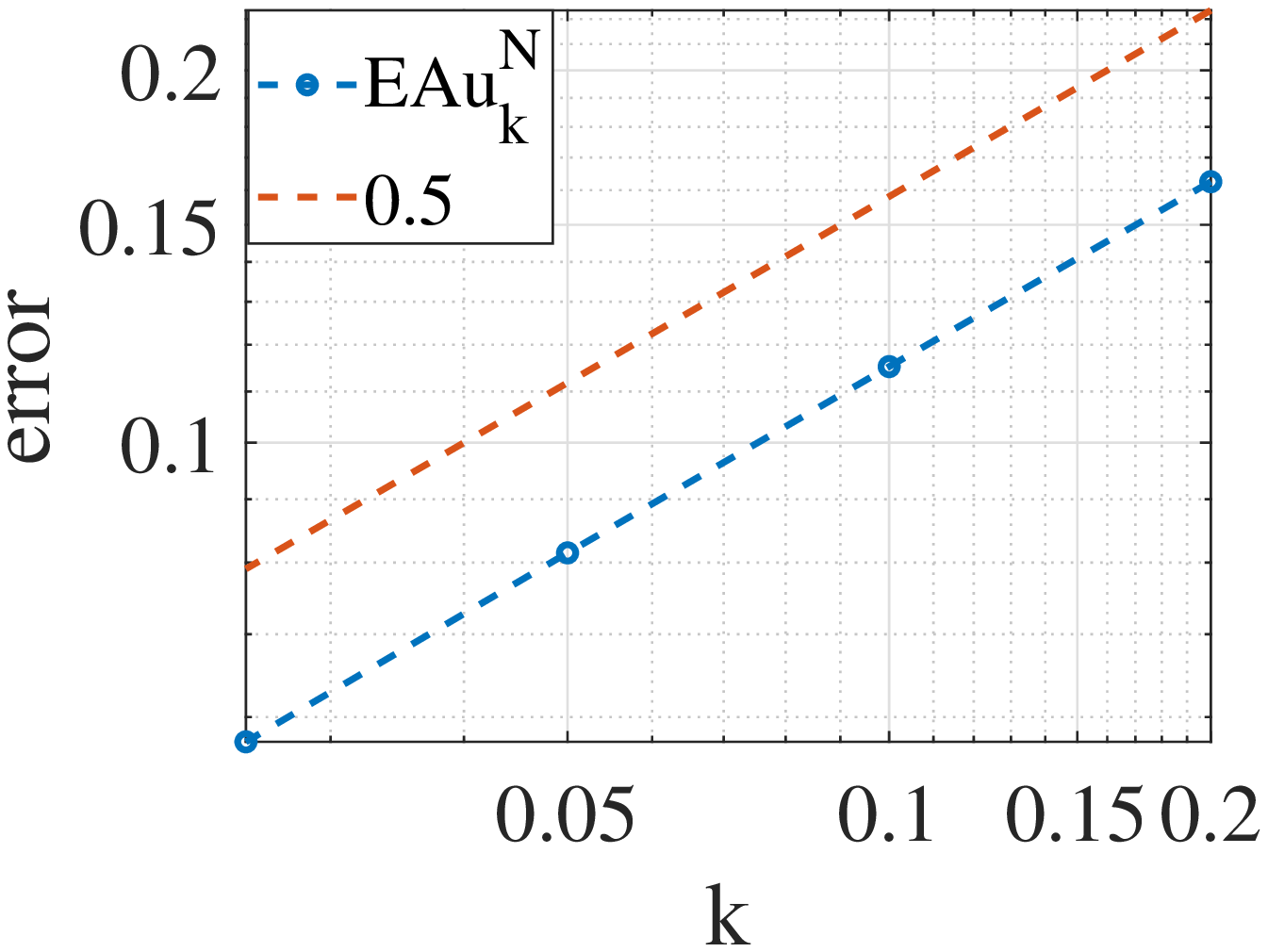}
}
\caption{Time errors of the numerical results at $T = 1$ with $h = 2^{-7}$.}
\label{fig1}
\end{figure}
\begin{figure}[th]
\centerline{
\includegraphics[scale=0.28]{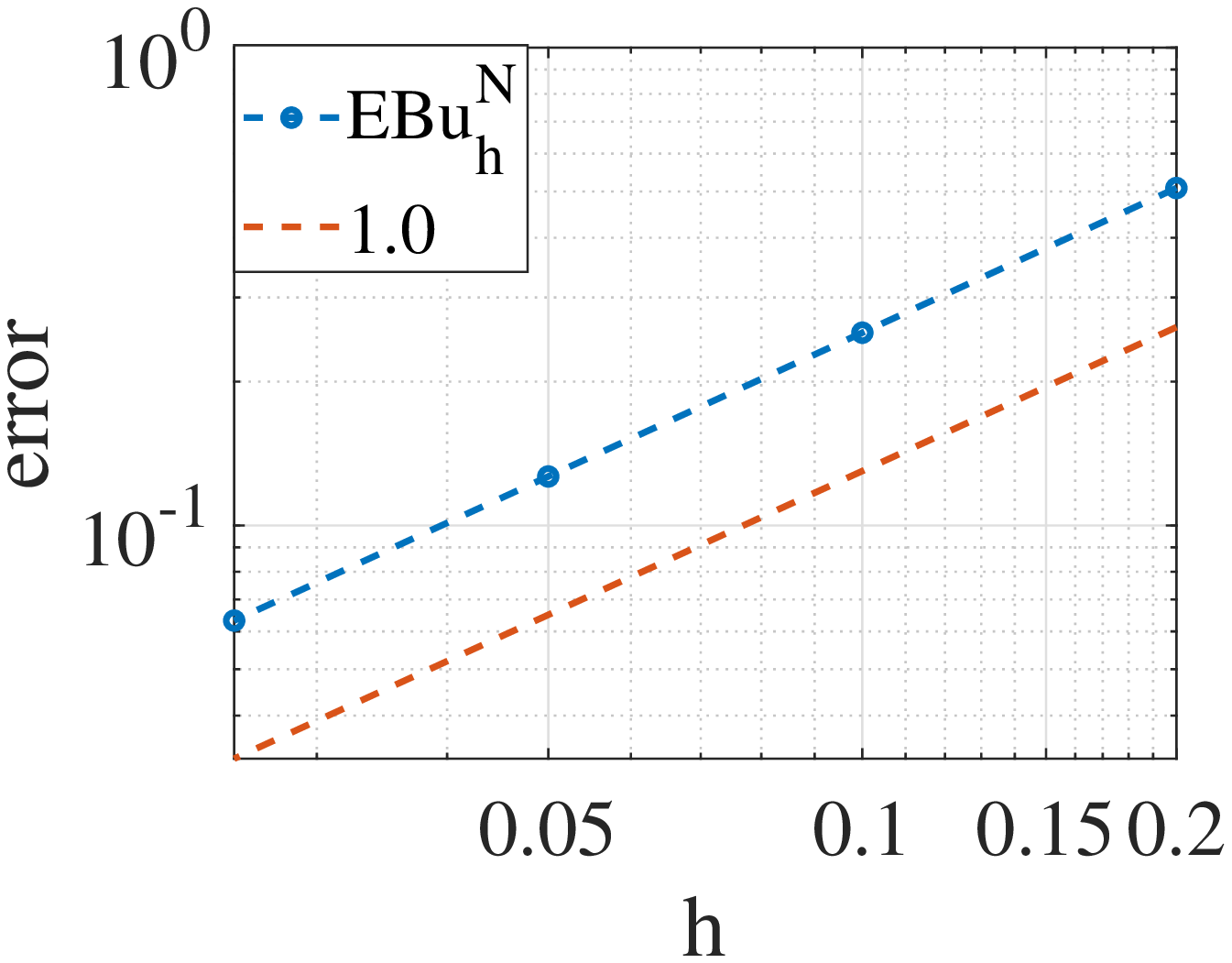}
\includegraphics[scale=0.28]{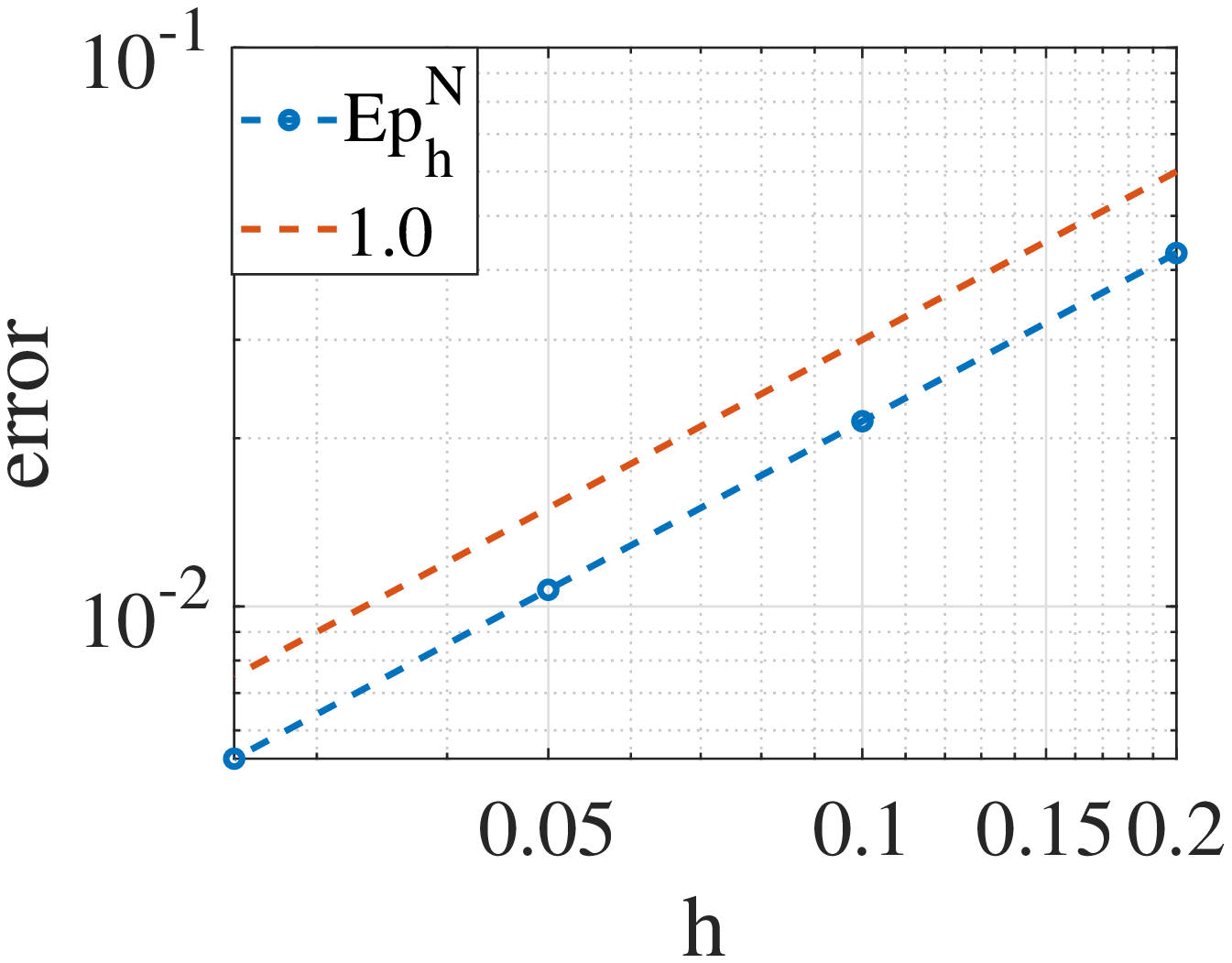}
\includegraphics[scale=0.28]{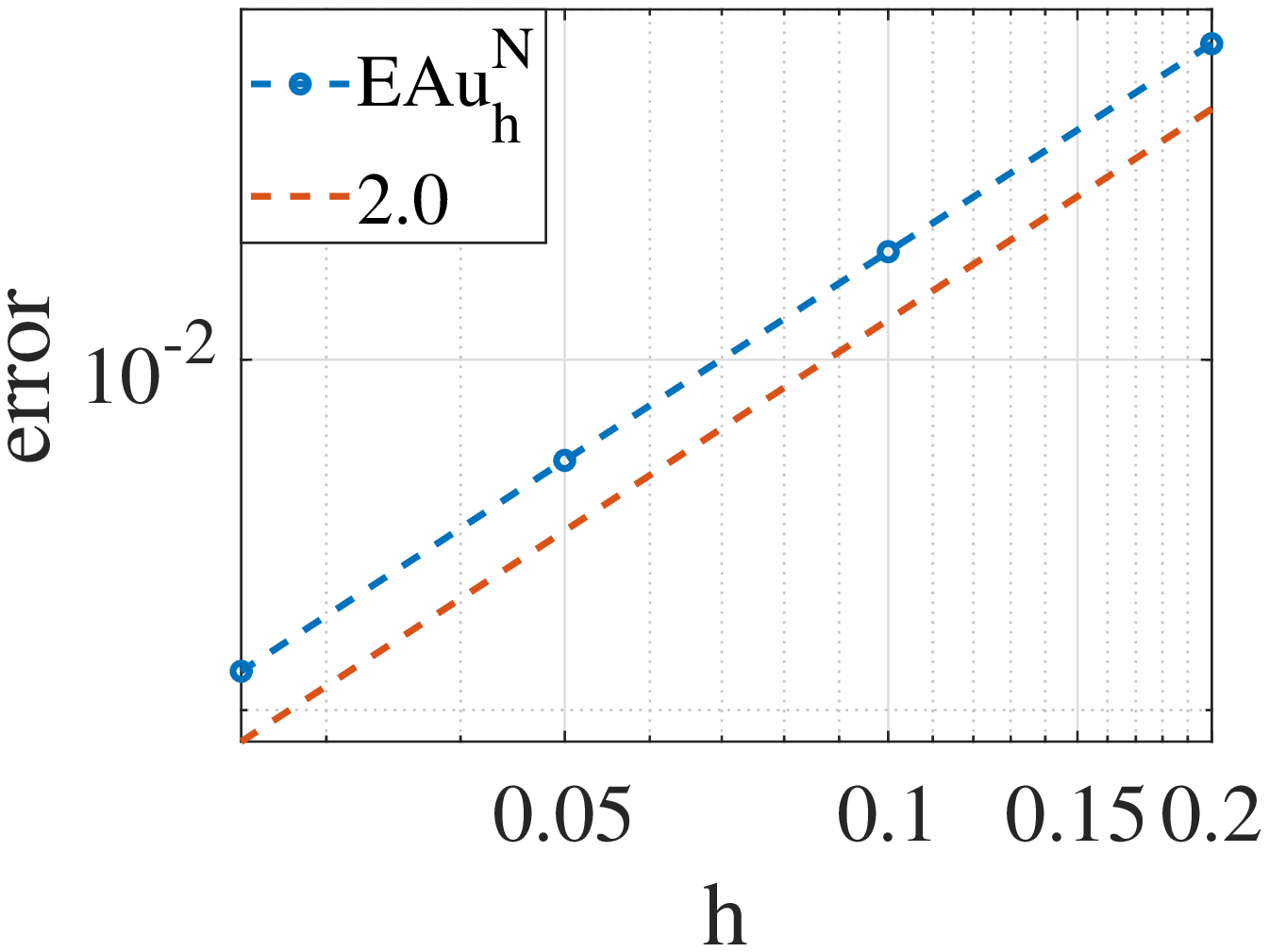}
}
\caption{Spatial errors of the numerical results at $T = 1$ with $k = 2^{-9}$.}
\label{fig2}
\end{figure}

\bigskip
\textbf{Acknowledgments.} The author would like to thank Professor Xiaobing Feng of The University of Tennessee
for his many discussions and critical comments as well as valuable suggestions which help to improve
the early version of the paper considerably. 

%
%



\begin{thebibliography}{99}


\bibitem{Bensoussan95} A. Bensoussan, {\em  Stochastic Navier-Stokes equations},
Acta Appl. Math., 38:267--304 (1995).

\bibitem{BT1973} A. Bensoussan and R. Temam, {\em Equations stochastiques du type Navier-Stokes},
J. Funct. Anal. 13, 195--222 (1973).


%

\bibitem{BBM14} H. Bessaih, Z. Brze\'{z}niak, and A. Millet, {\em Splitting up method for the 2D stochastic
	Navier-Stokes equations}, Stoch. PDE: Anal. Comp., 2:433--470 (2014).

\bibitem{BM18} H. Bessaih and A. Millet, {\em On strong $L^2$ convergence of time numerical schemes for the
	stochastic 2D Navier-Stokes equations}, IMA J. Numer. Anal., 39:2135--2167 (2019).

\bibitem{Breit21} D.  Breit and A. Dodgson, {\em Convergence rates for the numerical approximation of the
2D stochastic Navier-Stokes equations}, Numer. Math.,  147:553--578 (2021).

\bibitem{Brezzi_Fortin91} F. Brezzi and M. Fortin, {\em  Mixed and Hybrid Finite Element Methods},
Springer, New York, 1991.

\bibitem{BCP12} Z. Brze\'{z}niak, E. Carelli, and A. Prohl, {\em Finite element based discretizations of
the incompressible Navier-Stokes equations with multiplicative random forcing}, IMA J. Numer.
Anal., 33:771--824 (2013).

\bibitem{Capinski93} M. Capi\'{n}ski, {\em A note on uniqueness of stochastic Navier-Stokes equations}, Univ. Iagell. Acta Math.
30: 219--228 (1993).

\bibitem{Capinski91} M. Capi\'{n}ski,  N. J. Cutland, {\em Stochastic Navier-Stokes equations},  Acta Appl. Math. 25: 59--85 (1991).




\bibitem{CP12} E. Carelli and A. Prohl, {\em Rates of convergence for discretizations of the
stochastic incompressible Navier-Stokes equations}, SIAM J. Numer. Anal., 50(5):2467--2496 (2012).




\bibitem{Ern_Guermond04} A. Ern and J.-L. Guermond, {\em Theory and Practice of Finite Elements},
Springer, 2004.

\bibitem{D12} P. D\"{o}rsek, {\em Semigroup splitting and cubature approximations for
the stochastic Navier-Stokes equations}, SIAM J. Numer. Anal., 50(2):729--746 (2012).

\bibitem{Falk08} R. Falk, {\em A Fortin operator for two-dimensional Taylor-Hood elements},
ESAIM: Math. Model. Num. Anal., 42:411--424 (2008).

\bibitem{Feng_Qiu18} X. Feng and  H. Qiu, {\em Fully discrete mixed finite element methods for
the time-dependent stochastic Stokes equations with multiplicative noise},
J. Sci. Comput., 88:1-31 (2021).

\bibitem{Feng20} X. Feng, A. Prohl and L. Vo. {\em
Optimally convergent mixed finite element methods for the stochastic Stokes equations},
IMA J. Numer. Anal., 2021, doi.org/10.1093/imanum/drab006.


\bibitem{Flandoli_Gatarek95} F. Flandoli and D. Gatarek, {\em Martingale and stationary solutions for stochastic
Navier-Stokes equations}, Probab. Theory Related Fields, 102:367--391 (1995).

\bibitem{Hausenblas19} E. Hausenblas and T. Randrianasolo, {\em Time-discretization of stochastic
 2D Navier-Stokes equations with a penalty-projection method},
Numer. Math.,  143:339--378 (2019).




\bibitem{Girault_Raviart86}
V. Girault and P.-A. Raviart, {\em Finite Element Methods for Navier-Stokes Equations}, Springer, New York, 1986.

%
\bibitem{LRJ03} J. A. Langa, J. Real, and J. Simon, {\em Existence and regularity of the pressure for the stochastic
Navier-Stokes equations}, Appl. Math. Optim., 48:195--210 (2003).


\bibitem{H06}M. Hairer, J. C. Mattingly, {\em Ergodicity of the 2D Navier-Stokes equations with degenerate stochastic
forcing}, Ann. Math. 164, 993--1032 (2006).


\bibitem{K2012} S. Kuksin,  A. Shirikyan, {\em Mathematics of Two-Dimensional Turbulence}, Volume 194 of Cambridge
Tracts in Mathematics. Cambridge University Press, Cambridge, 2012.

\bibitem{PZ92} G. Da Prato and J. Zabczyk, {\em Stochastic Equations in Infinite Dimensions},
Cambridge University Press, Cambridge, UK, 1992.


\bibitem{Temam01} R. Temam,
{\em Navier-Stokes Equations. Theory and Numerical Analysis}, 2nd ed., AMS Chelsea
Publishing, Providence, RI, 2001.



\end{thebibliography}
\end{document}